\documentclass[a4paper,leqno]{amsart}
\usepackage[margin=30truemm]{geometry}
\usepackage{fancyhdr,appendix,graphicx,layout}
\usepackage{amsmath,amsthm,amssymb}
\usepackage[shortlabels]{enumitem}
\setlist[enumerate]{font=\normalfont}
\usepackage{physics}

\usepackage[numbers]{natbib}
\usepackage{bbm}
\usepackage{bm}
\usepackage{tikz}
\usepackage{tikz-cd}
\usetikzlibrary{positioning,arrows.meta,decorations.markings,calc}
\tikzset{
    cells={font=\everymath\expandafter{\the\everymath\displaystyle}},
}
\usepackage{caption}
\usepackage[pagewise]{lineno}

\usepackage[normalem]{ulem}
\usepackage{float}

\usepackage{color}
\definecolor{darkgreen}{rgb}{0,0.7,0}

\definecolor{darkblue}{rgb}{0,0,0.7} 
\newcommand{\defn}[1]{\textsl{\color{darkblue} #1}}

\usepackage[hyperpageref]{backref}
\usepackage[pagebackref]{hyperref}
\hypersetup{
    colorlinks=true,
    linkcolor=red,  
    filecolor=blue,  
    urlcolor=orange, 
    citecolor=violet   
}
\renewcommand*{\backref}[1]{}  
   \renewcommand*{\backrefalt}[4]{
      \ifcase #1 
         Not cited.
      \or
         Cited on page #2.
      \else
         Cited on pages #2.
      \fi}

\usepackage{cleveref}

\makeatletter
\newcommand\ReDeclareMathOperator[2]{%
    \begingroup \escapechar\m@ne\xdef\@gtempa{{\string#1}}\endgroup
    \expandafter\@ifundefined\@gtempa
        {\@latex@error{Command \string#1 undefined}\@ehc}
        \relax
    \let\@ifdefinable\@rc@ifdefinable
    \DeclareMathOperator#1{#2}}
\makeatother

\counterwithin{equation}{subsection}  
\crefname{thm}{Theorem}{Theorems}
\crefname{dfn}{Definition}{Definitions}
\crefname{dfnprop}{Definition-Proposition}{Definition-Proposition}
\crefname{prop}{Proposition}{Propositions}
\crefname{lem}{Lemma}{Lemmas}
\crefname{cor}{Corollary}{Corollaries}
\crefname{clm}{Claim}{Claims}
\crefname{ass}{Assumption}{Assumption}
\crefname{cond}{Condition}{Condition}
\crefname{nota}{Notation}{Notation}
\crefname{conj}{Conjecture}{Conjecture}
\crefname{fct}{Fact}{Facts}
\crefname{rmk}{Remark}{Remarks}
\crefname{eg}{Example}{Examples}
\crefname{que}{Question}{Questions}
\crefname{figure}{Figure}{Figures}
\crefname{table}{Table}{Tables}
\crefname{section}{Section}{Sections}
\crefname{subsection}{Subsection}{Subsections}
\crefname{appendix}{Appendix}{Appendices}
\crefname{equation}{}{}
\crefname{main}{Theorem}{Theorems}

\theoremstyle{plain}
\newtheorem{thm}{Theorem}[section]

\newtheorem{prop}[thm]{Proposition}
\newtheorem{lem}[thm]{Lemma}
\newtheorem{cor}[thm]{Corollary}

\theoremstyle{definition}
\newtheorem{dfn}[thm]{Definition}

\newtheorem{eg}[thm]{Example}

\theoremstyle{remark}
\newtheorem{rmk}[thm]{Remark}

\newcommand{\ZZ}{\mathbb{Z}}

\ReDeclareMathOperator{\top}{top}

\let\mod\relax
\DeclareMathOperator{\mod}{mod}

\DeclareMathOperator{\grmod}{grmod}

\DeclareMathOperator{\sdim}{\st{dim}}

\DeclareMathOperator{\Hom}{Hom}

\DeclareMathOperator{\Ext}{Ext}
\DeclareMathOperator{\Tor}{Tor}

\ReDeclareMathOperator{\l}{\ell}

\newcommand{\st}[1]{\underline{#1}}

\newcommand{\val}{\operatorname{val}}
\newcommand{\Max}{\operatorname{Max}}
\newcommand{\Ker}{\operatorname{Ker}}

\newcommand{\BOneN}{%
\left(\begin{smallmatrix}
I_n&0\\
0&0\\
0&0
\end{smallmatrix}\right)}

\newcommand{\AOneN}{%
\left(\begin{smallmatrix}
0&0&0\\
0&0&I_n
\end{smallmatrix}\right)}

\newcommand{\BTwoN}{%
\left(\begin{smallmatrix}
0&0\\
0&0\\
I_n&0
\end{smallmatrix}\right)}

\newcommand{\ATwoN}{%
\left(\begin{smallmatrix}
0&0&0\\
I_n&0&0
\end{smallmatrix}\right)}

\newcommand{\BThreeN}{%
\left(\begin{smallmatrix}
I_n&0\\
I_n&0\\
I_n&0
\end{smallmatrix}\right)}

\newcommand{\AThreeN}{%
\left(\begin{smallmatrix}
0&0&0\\
-J_n&I_n&J_n-I_n
\end{smallmatrix}\right)}

\DeclareMathOperator{\domdim}{domdim}

\DeclareMathOperator{\pdim}{pdim}
\DeclareMathOperator{\gldim}{gldim}

\DeclareMathOperator{\diam}{diam}

\DeclareMathOperator{\MH}{MH}

\renewcommand{\leq}{\leqslant}
\renewcommand{\le}{\leqslant}

\title[On distance algebras of finite graphs]{On representation-theoretic and homological properties of distance algebras of finite graphs}
\author{Aaron Chan and Bohan Xing}
\date{\today}

\newcommand{\Addresses}{{
  \bigskip
  \footnotesize

  A. Chan, \textsc{Department of mathematics, Nagoya University, Chikusa-ku, Nagoya 464-8602, Japan}\par\nopagebreak
  \textit{E-mail address}: \texttt{aaron.kychan@gmail.com}

  \medskip

  B. Xing, \textsc{School of Mathematical Sciences, Laboratory of Mathematics and Complex Systems, Beijing Normal University, Beijing 100875, P.R. China}\par\nopagebreak
\textit{E-mail address}: \texttt{bhxing@mail.bnu.edu.cn}
}}

\begin{document}

\begin{abstract}
Distance algebras were introduced by Asao and Ivanov to realize magnitude homology and cohomology as bigraded Tor and Ext groups.  We study the representation-theoretic and homological properties of the distance algebra $\Sigma_G$ of a finite connected graph $G$.  We characterize the graphs for which $\Sigma_G$ is representation-finite, multiserial, biserial, or self-injective by an explicit graph-theoretic classification.  We also compute the finitistic, self-injective, and dominant dimensions for several natural families of graphs.  For a finite modular lattice $L$, we consider the distance algebra $\Sigma_L$ associated to the Hasse graph of $L$, and show that $\Sigma_L$ is Auslander--Gorenstein precisely when $L$ is a finite direct product of chains.  We also present a family of minimal Auslander--Gorenstein distance algebras, namely, those associated to punctured hypercubes, whose finitistic (equivalently, self-injective; equivalently, dominant) dimension can be arbitrarily large.
\end{abstract}

\maketitle

\tableofcontents

\section{Introduction}

Distance algebras were introduced by Asao and Ivanov in \cite{AI24} to realize the magnitude (co)homology of finite quasi-metric spaces as bigraded Tor and Ext groups between the simple modules.  
A prominent class of examples of finite quasi-metric spaces is given by finite graphs, where the metric is given by the length of a shortest path.  This provides a natural link between magnitude (co)homology, graph geometry, and the representation theory of finite-dimensional algebras.  

As far as we can find in the literature, distance algebras have appeared in \cite{GHZ98} (without any special name), due to their connection with a class of non-commutative rings called \emph{tiled orders}; see \cite{GHZ98,IL26}  or Subsection \ref{subsec:tiled-order} for details.
The connection of these algebras with quasi-metric spaces was not utilized, and the study seems to not have gained mainstream attention.  
Motivated by the connection to magnitude (co)homology, we want to get a better understanding of this class of algebras, and this note presents some of the first investigations we have carried out. 

More precisely, for a finite connected graph $G$, let $\Sigma_G$ be the associated distance algebra (see Definition \ref{def:distance-algebra}).  
In general, $\Sigma_G$ could exhibit a wide range of representation-theoretic and homological properties.
Our aim is to determine the precise translation between several such properties and the combinatorics of $G$.

We also translate some standard techniques for computing magnitude (co)homology into statements about distance algebras 
in Section~\ref{sec:3-constr}.  Namely, Cartesian products of graphs correspond to tensor products of algebras, which is the algebraic counterpart of the K\"unneth theorem for magnitude homology~\cite{HW17}; convex full subgraphs $H\subseteq G$ correspond to corner algebras $e_H\Sigma_Ge_H$, and the projecting condition of Hepworth and Willerton~\cite{HW17} becomes a condition on modules over such a corner algebra; and we give an algebraic description of the distance algebra of a join of graphs (a construction considered in~\cite{HW17}).  While only some of them serve as tools in later sections, we feel that clarifying these translations between representation theory and topology will be beneficial to both communities for future research.

The first part of our main investigations concerns representation type and (multi-)seriality of $\Sigma_G$.  We write $A_n$ for the path graph  on $n$ vertices and $C_n$ for the cycle on $n$ vertices. 

\begin{thm}\textnormal{(see Theorem~\ref{thm:rep-finite} and
Propositions~\ref{prop:multiserial}, \ref{prop:biserial}, and
\ref{prop:biserial rep type})}
Let $G=(V,E)$ be a finite connected graph.  Then the following statements
hold.
\begin{enumerate}
    \item The algebra $\Sigma_G$ is representation-finite if and only if $G=A_n$ for some $n\geq1$.

    \item The following conditions are equivalent:
    \begin{enumerate}
        \item $\Sigma_G$ is multiserial;
        \item $\Sigma_G$ is special multiserial;
        \item $G$ is a path graph, a cycle, or a complete graph with a
        matching deleted.
    \end{enumerate}

    \item The following conditions are equivalent:
    \begin{enumerate}
        \item $\Sigma_G$ is biserial;
        \item $\Sigma_G$ is special biserial;
        \item $G=A_n$ for some $n\geq1$, or $G=C_n$ for some $n\geq3$.
    \end{enumerate}
\end{enumerate}
Moreover, $\Sigma_{C_3}$ and $\Sigma_{C_4}$ are representation-domestic, whereas $\Sigma_{C_n}$ is tame but not of polynomial growth for every $n\geq5$.
\end{thm}

Next, we will focus on homological properties of $\Sigma_G$.
Note that, except when $G$ consists of a single vertex, the associated distance algebra $\Sigma_G$ has infinite global dimension; see Proposition \ref{prop:cartan gldim}.
We first look at the graph-theoretic translation of the self-injectivity of $\Sigma_G$.  
We use the terminology of Section~\ref{sec:self-inj}, where $\Lambda_G$ denotes the tiled order associated with $G$ over the ring $R=\Bbbk[[x]]$ of formal power series.

\begin{thm}\textnormal{(see Theorem~\ref{thm:self-inj} and
Corollary~\ref{thm:tiled-gorenstein-antipodal})}
For a finite connected graph $G$, the following conditions are
equivalent:
\begin{enumerate}
    \item $G$ is interval-antipodal;
    \item $\Lambda_G$ is a 1-Iwanaga--Gorenstein $R$-order;
    \item $\Sigma_G$ is self-injective.
\end{enumerate}
\end{thm}

For finite modular lattices, the preceding criterion specializes to a particularly simple class: by Corollary~\ref{cor:modular-lattice-interval-antipodal}, the distance algebra $\Sigma_L$ of (the Hasse graph of) a finite modular lattice $L$ is self-injective if and only if $L$ is a Boolean lattice.  
In fact, we can also obtain a classification of Auslander--Gorenstein algebras (see Subsection \ref{subsec:def-homo-dim} for the definition) within such a class of distance algebras.

\begin{thm}\textnormal{(see Theorem~\ref{thm:modular-lattice-AG},
Corollary~\ref{cor:distributive-lattice-AG}, and Proposition~\ref{prop:distributive-lattice-IG})}
Let $L$ be a finite modular lattice.  The distance algebra $\Sigma_L$ is Iwanaga--Gorenstein if and only if $L$ is distributive. Moreover, the following conditions are equivalent.
\begin{enumerate}
    \item $\Sigma_L$ is Auslander--Gorenstein;
    \item $\operatorname{domdim}\Sigma_L\geq1$;
    \item $L$ is isomorphic to a finite direct product of chains.
\end{enumerate}
Moreover, when $L$ is a finite distributive lattice (hence, given by the lattice of order ideals of a finite poset $P$), the condition (3) can be replaced by
\begin{itemize}
    \item[\rm(3')] every connected component of $P$ is a chain.
\end{itemize}
\end{thm}

It turns out that to nail down the class of Iwanaga--Gorenstein distance algebras is generally a much harder problem than determining the Auslander--Gorenstein ones; in a companion note \cite{CX2026}, we will show that $\Sigma_L$ is Iwanaga--Gorenstein when $L$ is a finite distributive lattice.

\medskip 

We also compute the finitistic, self-injective, and dominant dimensions of the distance algebras associated with several natural families of graphs.  In Table \ref{tab:intro-homological-dimensions}, $T_n$ denotes an arbitrary tree on $n$ vertices, $K_n$ and $K_{m,n}$ denote the complete and complete bipartite graphs, respectively, $\square$ denotes the Cartesian product, and $Q_d^-$ denotes the punctured $d$-dimensional hypercube.

\begin{table}[htbp]
\centering
\caption{Homological dimensions of distance algebras for selected graph families.}
\label{tab:intro-homological-dimensions}
\small
\setlength{\tabcolsep}{4pt}
\renewcommand{\arraystretch}{1.15}
\begin{tabular}{c|c|c|c|c}
$G$ & condition & $\operatorname{findim}\Sigma_G$
& $\operatorname{idim}\Sigma_G$
& $\operatorname{domdim}\Sigma_G$ \\ \hline
$T_n$ & $n\leq2$ & $0$ & $0$ & $\infty$ \\
$T_n$ & $T_n=A_n$, $n\geq3$ & $1$ & $1$ & $1$ \\
$T_n$ & $T_n\neq A_n$, $n\geq3$ & $1$ & $1$ & $0$ \\ \hline
$C_n$ & $n\geq3$ even & $0$ & $0$ & $\infty$ \\
$C_n$ & $n\geq3$ odd & $0$ & $\infty$ & $0$ \\ \hline
$K_n$ & $n\geq3$ & $0$ & $\infty$ & $0$ \\ \hline
$K_{m,n}$ & $m,n\geq2$ and $(m,n)\neq(2,2)$
& $0$ & $\infty$ & $0$ \\ \hline
$K_{2\times n}$ & $n\geq2$
& $0$ & $0$ & $\infty$ \\ \hline
$A_{k_1}\square\cdots\square A_{k_n}$
& $k_i\geq3$ for $1\leq i\leq n$ & $n$ & $n$ & $1$ \\ \hline
$Q_d^-$ & $d\geq2$ & $d-1$ & $d-1$ & $d-1$
\end{tabular}
\end{table}

In particular, $\Sigma_{Q_d^-}$ is $(d-2)$-minimal Auslander--Gorenstein for every $d\geq2$.  We also describe, under the Auslander--Iyama--Solberg correspondence, the associated precluster tilting module in Subsection~\ref{subsubsec:precluster}.

\medskip
\noindent
\textbf{Organization of the paper.}
Section~\ref{sec:2} recalls the basic properties of distance algebras, the homological dimensions used throughout, and graph-metric tiled orders.  Section~\ref{sec:3-constr} develops elementary constructions relating graphs and their distance algebras.  Section~\ref{sec:4-proj-dim} studies the projective dimensions of path-generated right ideals.  Section~\ref{sec:biserial} treats representation-finiteness,
multiseriality, and biseriality.  Section~\ref{sec:self-inj} studies
self-injectivity, tiled orders, and the relevant lattice-theoretic
constructions.  Finally, Section~\ref{sec:tree} gives the homological dimension calculations and the applications to modular and distributive lattices and punctured hypercubes.

\medskip
\noindent
{\bf Conventions and notation.}
For a path $p$ in a quiver, we denote by $s(p)$ its source vertex and by $t(p)$ its target vertex, and by $\ell(p)$ its length. We compose arrows from left to right, and map compositions are read from right to left. All modules are right modules.

\medskip
\noindent
\textbf{Use of AI.}
ChatGPT model 5.5 pointed out Example \ref{eg:self-centered eg} can be unified using self-centered graphs.  It also suggested the punctured hypercube as a candidate minimal Auslander--Gorenstein algebra to us.  We have used ChatGPT model 5.6 Sol for generalizing our results and arguments from distributive lattices to modular lattices.
All mathematical arguments produced by AI were independently verified and finalized by the authors. 
ChatGPT has helped to improve the writing throughout the paper.  
The authors take full responsibility for the accuracy, originality, and integrity of the paper.

\medskip
\noindent
\textbf{Acknowledgments.}
We thank Osamu Iyama for telling us about the connection between distance algebras and tiled orders.  AC is deeply grateful to Yasuhiko Asao and Shun Wakatsuki for teaching him about the connection between distance algebras and magnitude homologies, and the many discussions on the subject.
We thank Ren\'{e} Marczinzik for pointing out Proposition~\ref{prop:distributive-lattice-IG}(3).
AC is supported by JSPS Grant-in-Aid for Scientific Research (C) 24K06666. BX is supported by the China Scholarship Council (No. 202506040127).

\section{Preliminaries}\label{sec:2}

\subsection{Distance algebras of graphs} 

Let $G=(V,E)$ be a connected combinatorial graph, where $V$ is the set of vertices and $E$ is the set of edges. In this paper, we assume that $G$ contains neither loops nor multiple edges.

For vertices $x,y\in V$, we write $\{x,y\}\in E$ if $x$ and $y$ are connected by an edge; in this case, they are called \defn{adjacent}. The \defn{distance} between $x$ and $y$, denoted by $$d(x,y)=d_G(x,y),$$ is the length of a shortest path from $x$ to $y$. Recall that 
\begin{itemize}
    \item the \defn{valency} of a vertex $x$ is the number of vertices adjacent to $x$, denoted by $$\val(x)=\operatorname{val}_G(x):=|\{y\in V\mid \{x,y\}\in E\}|;$$

    \item the \defn{eccentricity} of a vertex $x$ is the largest distance from $x$, denoted by $$\operatorname{ecc}(x):=\max_{y\in V} d(x,y);$$

    \item the \defn{diameter} of $G$ is the maximal distance in $G$, or equivalently, the maximal eccentricity, denoted by $$\operatorname{diam}(G) :=\max_{x,y\in V}d(x,y)=\max_{x\in V}\operatorname{ecc}(x).$$
\end{itemize}

\begin{dfn}\label{def:distance-algebra}
    Let $G=(V,E)$ be a graph (without multiple edges and loops), and $Q_G$ be its \defn{double quiver}, i.e., the set of vertices of $Q_G$ is $V$ and there is a pair of arrows 
\[\begin{tikzcd}
	x \arrow[r, "{a_{xy}}", bend left] & y \arrow[l, "{a_{yx}}", bend left]
\end{tikzcd} \text{ for each edge }\{x,y\}\in E.\]
Let $\Bbbk$ be a unital commutative ring.
The \defn{distance $\Bbbk$-algebra} of $G$ is defined to be
$\Sigma_G:=\Bbbk Q_G/I_G$, where $I_G$ is the ideal generated by the following relations:
\begin{enumerate}
    \item $p=0$, for every path $p$ such that $\ell(p)>d(s(p),t(p))$;
    \item $p=q$, for every pair of paths $p,q$ in $Q_G$ with
$$s(p)=s(q),\qquad t(p)=t(q),\qquad
\ell(p)=d(s(p),t(p))=\ell(q).$$
\end{enumerate}
\end{dfn}

By definition, in $\Sigma_G$, a path is non-zero if and only if it is represented by a shortest path in $G$, and any two shortest paths with the same starting and ending vertices are identified. Hence, for any two vertices $x,y\in V$, there is a unique non-zero path from $x$ to $y$ modulo $I_G$. 
In particular, $\Sigma_G$ is Schurian in the sense that 
$\dim_{\Bbbk} e_x\Sigma_Ge_y=1$ for all $x,y\in V$, and so $\dim_{\Bbbk} \Sigma_G=|V|^2$.

We take a $\Bbbk$-basis $\{b_{xy}\mid x,y\in V\}$ for $\Sigma_G$,
where $b_{xy}$ denotes the equivalence class of any shortest path from $x$ to $y$ on $G$, with the convention that $b_{xx}=e_x$ is the primitive idempotent associated to $x$. Then multiplication in $\Sigma_G$ is given by
$$
b_{xy}b_{yz}=
\begin{cases}
b_{xz}, & \text{if } d(x,z)=d(x,y)+d(y,z),\\
0, & \text{otherwise}.
\end{cases}
$$
Note that $b_{xy}$ is denoted by $(x,y)$ and by $e_{xy}$ in \cite{AI24} and \cite{AW24}, respectively. 

\medskip

Note that path length induces a natural grading on $\Sigma_G$; this matches the canonical (internal $\mathbb{R}$-)grading of magnitude homology.  
The following observation is straightforward.

\begin{lem}\label{lem:Loewy legnth=ecc}
Denote by $\operatorname{LL}(M)$ the Loewy length of a $\Sigma_G$-module $M$.  We have the following equalities.
\begin{enumerate}
    \item $\operatorname{LL}(P_x)=\operatorname{ecc}(x)+1$.
    \item $\operatorname{LL}(\Sigma_G)=\operatorname{diam}(G)+1$. 
\end{enumerate}
\end{lem}

This class of algebras recently appeared in a study by Asao and Ivanov on categorifying magnitude (co)homology \cite[Theorem 6.2]{AI24}; they show that the (graded) Tor groups and Ext groups of the $\Sigma_G$-module $S_G:=\bigoplus_{v\in V(G)} S_v$ coincide with the $\Bbbk$-valued magnitude homology $\MH_\bullet^*(G;\Bbbk)$ and cohomology $\MH^{\bullet,*}(G;\Bbbk)$, respectively, for any unital commutative ring $\Bbbk$:
\begin{align*}
\MH_n^
\ell(G;\Bbbk) & \cong \Tor_{n,\ell}^{\Sigma_G}(S_G,S_G):=\Tor^{\grmod\Sigma_G}_{n}(S_G,S_G(\ell)), \\
\MH_n(G;\Bbbk) & \cong \Tor^{\Sigma_G}_{n}(S_G,S_G), \\
\MH^{n,\ell}(G;\Bbbk) & \cong \Ext^{n,\ell}_{\Sigma_G}(S_G,S_G):=\Ext_{\grmod\Sigma_G}^{n}(S_G,S_G(\ell)), \\
\MH^{n}(G;\Bbbk) & \cong \Ext_{\Sigma_G}^{n}(S_G,S_G).
\end{align*}
Here $(1)$ denotes the grading shift functor $M(1)_i := M_{i-1}$. 

A graph $G$ is \defn{diagonal} if $\MH_{n,\ell}(G;\ZZ)=0$ whenever $n\neq \ell$.
By the K\"{u}nneth theorem, this is equivalent to $\operatorname{MH}_{n,\ell}(G;\Bbbk)=0$ whenever $n\neq\ell$ over all fields $\Bbbk$.
As a consequence of the above identification between Tor-groups and magnitude homologies, $\MH_{n,\ell}(G;\Bbbk)$ being concentrated in the diagonal is equivalent to saying that the $\Bbbk$-algebra $\Sigma_G$ is a so-called \emph{Koszul algebra}.  We refer to \cite{IM24} for the details.

For graph-theoretic notation and terminology, we follow
\cite{BM08,HIK11}. The main families of graphs considered in this paper are listed below, together with references to the corresponding results on their distance algebras.
\begin{itemize}
    \item $T_n$ denotes an arbitrary tree with $n$ vertices.
    See Subsection~\ref{subsec:tree} for the homological dimensions
    of $\Sigma_{T_n}$.

    \item $A_n$ denotes the path graph with $n$ vertices.
    The algebra $\Sigma_{A_n}$ is representation-finite by
    Theorem~\ref{thm:rep-finite}.
    Its homological dimensions are computed in
    Subsection~\ref{subsec:tree}.

    \item $C_n$ denotes the cycle with $n\geq3$ vertices.
    See Proposition~\ref{prop:biserial rep type} for the tame
    representation type of $\Sigma_{C_n}$.
    For even $n$, the algebra $\Sigma_{C_n}$ is self-injective;
    see Example~\ref{ex:interval-antipodals}(1).
    For odd $n$, its homological dimensions are computed in
    Subsection~\ref{sec:self-centered}.

    \item $K_n$ denotes the complete graph with $n$ vertices.
    See Subsection~\ref{sec:self-centered} for the homological
    dimensions of $\Sigma_{K_n}$.

    \item $K_{m,n}$ denotes the complete bipartite graph with parts
    of cardinalities $m$ and $n$.
    See Subsection~\ref{sec:self-centered} for the homological
    dimensions of $\Sigma_{K_{m,n}}$ when $m,n\geq2$, and
    Subsection~\ref{subsec:tree} when $\min\{m,n\}=1$.
    The self-injective case $K_{2,2}\cong C_4$ is covered by     Example~\ref{ex:interval-antipodals}(1).

    \item
    $K_{2\times n}=
    K_{{2,\ldots,2}}$, for $n\geq2$,
    denotes the complete $n$-partite graph with two vertices in
    each part, also called the \defn{cocktail party graph}.
    Equivalently, it is obtained from $K_{2n}$ by deleting a
    perfect matching.
    Its distance algebra is self-injective; see
    Proposition~\ref{prop:diam2 int-antipodal <=> cocktail party}.

    \item $Q_n$ denotes the $n$-dimensional hypercube with vertex
    set $2^{[n]}$, where $[n]:=\{1,2,\ldots,n\}$.
    Two subsets $S,T\subseteq[n]$ are adjacent precisely when
    $|S\mathbin{\triangle}T|=|(S\setminus T)\cup(T\setminus S)|=1$.
    The algebra $\Sigma_{Q_n}$ is self-injective; see
    Example~\ref{ex:interval-antipodals}(2).

    \item $Q_n^-:=Q_n\setminus\{[n]\}$ denotes the
    \defn{punctured $n$-dimensional hypercube}, equipped with
    the metric inherited from $Q_n$.
    See Subsection~\ref{sec:ex-2-IG} for the homological
    dimensions of $\Sigma_{Q_n^-}$.

    \item $Q_n^\pm:=Q_n\setminus\{\emptyset,[n]\}$, for $n\geq3$,
    denotes the \defn{doubly punctured $n$-dimensional hypercube},
    equipped with the metric inherited from $Q_n$.
    The algebra $\Sigma_{Q_n^\pm}$ is self-injective; see
    Example~\ref{ex:interval-antipodals}(2).
\end{itemize}

\subsection{Homological dimensions}\label{subsec:def-homo-dim}

For a right $A$-module $M$, we denote by $\operatorname{pdim}(M)$ the projective dimension of $M$ and by $\operatorname{idim}(M)$ the injective dimension of $M$, respectively.

The \defn{global dimension} of $A$, i.e., the supremum of projective (equivalently, injective) dimensions across all $A$-module, is denoted by $\gldim A$.
The \defn{finitistic dimension} of $A$ is 
\[
\operatorname{findim}(A):=\sup \{ \pdim(M) \mid M \text{ has finite projective dimension} \}.
\]
The \defn{big finitistic dimension} of \(A\) is
$$\operatorname{Findim}(A)
:=\sup\{\pdim(M)\mid M \text{ is an arbitrary \(A\)-module of finite projective dimension}\}.$$
For a general algebra, one can distinguish the  finitistic projective dimension (the one above) and the finitistic injective dimension (equivalently, right and left finitistic dimension).  In the setting in which we are interested, it turns out that the two are the same (Corollary \ref{cor:fdim symmetry}).

An algebra is \defn{Iwanaga--Gorenstein} if both the injective dimension of $A$ as a left $A$-module and as a right $A$-module are finite.  In this case, it follows from Zaks' lemma that these two numbers coincide, and we call it the \defn{self-injective dimension} of $A$, and say that $A$ is $d$-Iwanaga--Gorenstein if the self-injective dimension is $d$.  Note that, in such a case, we have
$\operatorname{findim}(A)
=d=
\operatorname{Findim}(A)$
by \cite[Theorem~3.2]{AHHT}.

Let $0\to M\to I^0 \to I^1 \to \cdots$ be a minimal injective resolution of $M$.  Then we write $\domdim(M)\ge n$ if $I^i$ is projective for all $i<n$.  The \defn{dominant dimension} of $M$ is the supremum of all $n\in\mathbb{Z}_{\ge 0}$ such that $\domdim(M)\ge n$; in which case, we write $\domdim(M)=n\in \mathbb{Z}_{\ge 0}\cup \{\infty\}$.  Note that each projective-injective module has infinite dominant dimension.
The dominant dimension of an algebra $A$ is the dominant dimension of $A$ as the regular representation over itself, and this is the same number even if one works over the left modules.

Let $d\ge 0$ be a non-negative integer.
An algebra $A$ is called
\begin{itemize}
    \item \defn{Auslander--Gorenstein} if $A$ is Iwanaga--Gorenstein and $\pdim I^k\le k$ for all $k\ge 0$, where $0\to A\to I^0 \to I^1 \to I^2\to \cdots$ is the minimal injective coresolution of $A_A$.

    \item \defn{minimal Auslander--Gorenstein}, or more specifically \defn{$d$-minimal Auslander--Gorenstein}, if it is Iwanaga--Gorenstein and $\operatorname{idim} A \le d+1\le \domdim A$.
\end{itemize}
Note that if one assumes $A$ is ring-indecomposable and non-semisimple, then the two inequalities above can be replaced by equalities.

\subsection{Tiled orders}\label{subsec:tiled-order}

Throughout this subsection, unless otherwise stated, we assume that $\Bbbk$ is a field.
Let $R:=\Bbbk[[x]]$ be the formal power series ring and let $K:=\Bbbk((x))$ be its field of fractions.  Recall that a
\defn{tiled $R$-order} in $M_n(K)$ is an $R$-order of the form
\[
\Lambda=(x^{m_{ij}}R)_{1\leq i,j\leq n}\subseteq M_n(K),
\]
where $m_{ii}=0$ and $m_{ij}+m_{jk}\geq m_{ik}$ for all $i,j,k$.  These inequalities are precisely what is needed for $\Lambda$ to be closed under matrix multiplication.

For a finite connected graph $G=(V=\{v_1,\ldots, v_n\},E)$, its associated tiled order is
\[
\Lambda_G
:=
\begin{pmatrix}
R & x^{d(v_1,v_2)}R & \cdots & x^{d(v_1,v_n)}R\\
x^{d(v_2,v_1)}R & R & \cdots & x^{d(v_2,v_n)}R\\
\vdots & \vdots & \ddots & \vdots\\
x^{d(v_n,v_1)}R & x^{d(v_n,v_2)}R & \cdots & R
\end{pmatrix}
\subseteq M_n(K).
\]
The triangle inequality shows that $\Lambda_G$ is a tiled $R$-order.
Moreover, it is basic, since $d(v_i,v_j)+d(v_j,v_i)>0$ whenever $i\neq j$.

The following observation identifies the distance algebra as the
special fiber of $\Lambda_G$.

\begin{prop}\label{prop:graph-metric-tiled-order}
There are natural isomorphisms
\[
\Lambda_G/x\Lambda_G\cong\Sigma_G
\qquad\text{and}\qquad
\Lambda_G\otimes_R K\cong M_n(K).
\]
The first isomorphism preserves the natural gradings.
\end{prop}

\begin{proof}
Let $E_{i,j}$ be the matrix with $1$ at the $(i,j)$-entry and $0$ everywhere else.  The residue classes
$\overline{x^{d(v_i,v_j)}E_{ij}}$, for $1\leq i,j\leq n$,
form a $\Bbbk$-basis of $\Lambda_G/x\Lambda_G$.  Define
\[
\Phi:\Lambda_G/x\Lambda_G\longrightarrow\Sigma_G,
\qquad
\overline{x^{d(v_i,v_j)}E_{ij}}\longmapsto b_{v_i v_j}.
\]
For $u,v,w\in V$, one has
\[
(x^{d(u,v)}E_{uv})(x^{d(v,w)}E_{vw})
=x^{d(u,v)+d(v,w)-d(u,w)}
 (x^{d(u,w)}E_{uw}).
\]
The exponent outside the parentheses is non-negative by the triangle
inequality.  Consequently, the product vanishes modulo
$x\Lambda_G$ unless
$d(u,w)=d(u,v)+d(v,w)$, in which case it is sent by $\Phi$ to
$b_{uw}$.  This is exactly the multiplication rule in $\Sigma_G$,
so $\Phi$ is a graded algebra isomorphism.  The second isomorphism
follows because $x$ becomes invertible in $K$.
\end{proof}

The symmetry of the graph distance gives compatible dualities for the order and its special fiber.

\begin{lem}\label{lem:anti-iso}
The transpose map restricts to an $R$-algebra isomorphism
\[
\omega_R:\Lambda_G\longrightarrow\Lambda_G^{\operatorname{op}},
\qquad
x^{d(u,v)}E_{uv}\longmapsto x^{d(u,v)}E_{vu}.
\]
Reducing modulo x gives a $\Bbbk$-algebra isomorphism
\[
\omega_{\Bbbk}:\Sigma_G\longrightarrow
\Sigma_G^{\operatorname{op}},
\qquad
b_{uv}\longmapsto b_{vu}.
\]
The latter isomorphism holds over any unital commutative ring $\Bbbk$.
\end{lem}

\begin{proof}
Since $d(u,v)=d(v,u)$, transposition preserves $\Lambda_G$ and
reverses multiplication.
\end{proof}

As a consequence, $\omega$ induces  a simple-preserving duality $\omega_*:\mod \Sigma_G \to \mod \Sigma_G^{\operatorname{op}}$.

\begin{prop}\label{prop:dual-functor}
Let $\omega:\Sigma_G\rightarrow \Sigma_G^{\operatorname{op}}$
be the algebra isomorphism given above and let $D=\operatorname{Hom}_{\Bbbk}(-,\Bbbk)$ be the $\Bbbk$-linear duality. Then
$$\omega_*^{-1}D:\operatorname{mod} \Sigma_G\longrightarrow \operatorname{mod} \Sigma_G$$ is a duality. Moreover, for every finite-dimensional right $\Sigma_G$-module $M$, 
$$\operatorname{idim}M=\operatorname{pdim}(\omega^{-1}DM),
\qquad\operatorname{pdim}M=\operatorname{idim}(\omega^{-1}DM).$$
\end{prop}

\begin{proof}
It follows from the properties of the simple-preserving duality $\omega$ and the exact contravariant $\Bbbk$-duality $D=\operatorname{Hom}_{\Bbbk}(-,\Bbbk)$.
\end{proof}

The following is an immediate consequence.
\begin{cor}\label{cor:fdim symmetry}
    $\omega_*^{-1}D$ restricts to an equivalence between the full subcategory consisting of modules of finite projective dimension and that of modules of finite injective dimension.
\end{cor}

We next record the resulting relation between homological dimensions.
For $\Lambda_G$, the finitistic dimension is taken over finitely
generated right modules.  By Lemma \ref{lem:anti-iso}, the left and right versions of the dimensions below coincide.

\begin{prop}\label{prop:homological-dimensions-tiled-order}
Let $G$ be a finite connected graph.  Then
\[
\operatorname{idim}\Lambda_G
=\operatorname{idim}\Sigma_G+1, \qquad
\operatorname{findim}\Lambda_G
=\operatorname{findim}\Sigma_G+1.
\]
If $\Lambda_G$ has finite self-injective dimension  $\operatorname{idim}\Lambda_G<\infty$, then
$$\operatorname{Findim}\Lambda_G
=\operatorname{findim}\Lambda_G
=\operatorname{idim}\Lambda_G
=\operatorname{Findim}\Sigma_G+1
=\operatorname{findim}\Sigma_G+1
=\operatorname{idim}\Sigma_G+1.$$
In particular, if $\operatorname{gldim}\Lambda_G=d+1<\infty$ for some $d\ge 0$, then
\[
\operatorname{Findim}\Sigma_G=\operatorname{findim}\Sigma_G=\operatorname{idim}\Sigma_G=d.
\]
\end{prop}
This is a standard exercise in the homological algebra of orders; we give a brief explanation for completeness.
\begin{proof}
It suffices to consider right modules. Since $x$ is a central non-zerodivisor in $\Lambda_G$ and $\Lambda_G/x\Lambda_G\cong\Sigma_G\neq0$, the second injective change-of-rings theorem in \cite[Exercise~4.3.3]{Wei94} gives
\[
 \operatorname{idim}_{\Lambda_G}\Lambda_G
 \geq\operatorname{idim}_{\Sigma_G}\Sigma_G+1.
\]
Suppose $n=\operatorname{idim}_{\Sigma_G}\Sigma_G<\infty$. The first theorem in the same exercise gives $\operatorname{idim}_{\Lambda_G}\Sigma_G=n+1$. For any finitely generated right $\Lambda_G$-module $M$, the exact sequence
\[
 0\longrightarrow\Lambda_G\xrightarrow{\cdot x}\Lambda_G
 \longrightarrow\Sigma_G\longrightarrow0
\]
therefore yields
\[
 \operatorname{Ext}^i_{\Lambda_G}(M,\Lambda_G)
 =x\operatorname{Ext}^i_{\Lambda_G}(M,\Lambda_G)
 \qquad(i>n+1).
\]
These Ext groups are finitely generated over $\Bbbk[[x]]$, so they vanish by Nakayama's lemma. Since $\Lambda_G$ is noetherian, this implies $\operatorname{idim}_{\Lambda_G}\Lambda_G\leq n+1$, proving the first equality, including the infinite case.

The finitistic dimension formula is \cite[Proposition~4.1]{GHZ98}. When $\Lambda_G$ has finite self-injective dimension, the displayed chain of equalities follows from these two formulas and \cite[Propositions~1.2 and~1.3]{GHZ98}. Finally, if $\operatorname{gldim}\Lambda_G=d+1<\infty$, then $\operatorname{Findim}\Lambda_G=d+1$, giving the last assertion.
\end{proof}


\section{Elementary constructions in graphs versus in algebras}\label{sec:3-constr}

Unless otherwise stated, we keep our assumption that $\Bbbk$ is a unital commutative ring.

\subsection{Cartesian product versus tensor product}
\begin{dfn}
Let $G$ and $H$ be two graphs. The \defn{Cartesian product} of $G$ and $H$, denoted by $G\square H$, is the graph whose vertex set is
$$V(G\square H)=V(G)\times V(H).$$
Two vertices $(x,u)$ and $(y,v)$ are adjacent in $G\square H$ if and only if
\[
\text{either }\begin{cases}
x=y\\ \{u,v\}\in E(H),
\end{cases}\; \text{or }
\begin{cases}
 \{x,y\}\in E(G)\\u=v.
\end{cases}
\]
\end{dfn}

\begin{prop}\label{prop:tensor}
We have an algebra isomorphism
$\Sigma_{G\square H}\cong \Sigma_G\otimes_{\Bbbk} \Sigma_H$.
\end{prop}

\begin{proof}
This result can also be obtained directly from the construction of tensor products of quiver algebras in \cite[Proposition 3]{Her08}. We include here a concrete verification using the distance function.

For vertices $(x,u),(y,v)\in V(G)\times V(H)$, the distance in the Cartesian product is given by
$$d_{G\square H}((x,u),(y,v))=d_G(x,y)+d_H(u,v).$$
Now define a $\Bbbk$-linear map $\Phi:\Sigma_{G\square H}\longrightarrow \Sigma_G\otimes_{\Bbbk} \Sigma_H$ by linearly extending
$$b^{G\square H}_{(x,u),(y,v)}\mapsto b^G_{xy}\otimes b^H_{uv}.$$
Since both sides have bases indexed by pairs of vertices in $V(G)\times V(H)$, the map $\Phi$ is a $\Bbbk$-linear isomorphism.

It is straightforward to check that $b^{G\square H}_{(x,u),(y,v)}b^{G\square H}_{(y,v),(z,w)}$ is non-zero, with product being $b^{G\square H}_{(x,u),(z,w)}$, precisely when
$$d_G(x,z)+d_H(u,w)=d_G(x,y)+d_H(u,v)+d_G(y,z)+d_H(v,w).$$
By the triangle inequality in $G$ and $H$, this holds if and only if both $d_G(x,z)=d_G(x,y)+d_G(y,z)$ and $d_H(u,w)=d_H(u,v)+d_H(v,w)$ hold.
Since $(b^G_{xy}\otimes b^H_{uv})(b^G_{yz}\otimes b^H_{vw})=(b^G_{xy}b^G_{yz})\otimes (b^H_{uv}b^H_{vw})$ by the definition of the tensor product algebra, it follows from the multiplication rule that $\Phi(bb')=\Phi(b)\Phi(b')$ for all basis elements.
Thus, $\Phi$ is an algebra isomorphism.
\end{proof}


\begin{prop}
The enveloping algebra $\Sigma_G\otimes_\Bbbk \Sigma_G^{\operatorname{op}}$ of $\Sigma_G$ is isomorphic to the distance algebra of the Cartesian product of $G$ with itself.
\end{prop}
\begin{proof}
    Combine Lemma \ref{lem:anti-iso} and Proposition \ref{prop:tensor}.
\end{proof}

\subsection{Convex full subgraph versus idempotent subalgebra}\label{subsec:full subgraph}

For a subset $W\subseteq V$, the \defn{full subgraph} of $G$ induced by $W$, denoted by $G|_W$, is the graph with vertex set $W$ and edge set
$$E(G|_W)=E(G)\cap (W\times W)=\bigl\{\{x,y\}\in E(G)\mid x,y\in W\bigr\}.$$

We say that a subgraph $H$ of $G$ is \defn{convex}  if $d_H(x,y)=d_G(x,y)$ for all $x,y\in V(H)$.  Note that every non-empty convex subgraph of a connected graph is connected.

For a subgraph $H$ of $G$, we define an idempotent
\[e_H:=\sum_{x\in V(H)}e_x\in \Sigma_G.
\]

\begin{lem}\label{lem:corner-convex-subgraph}
Let \(G=(V,E)\) be a graph.  Let $H$ be the full subgraph of $G$ defined by  \(W\subseteq V\).
If $H$ is convex, then there is a natural algebra isomorphism
\[\Sigma_H\cong e_H\Sigma_Ge_H.\]
\end{lem}

\begin{proof}
For \(x,y\in W\), denote by \(b^H_{xy}\) the basis element of \(\Sigma_H\), and by \(b^G_{xy}\) the basis element of \(\Sigma_G\). Since $e_H\Sigma_Ge_H=\operatorname{span}_{\Bbbk}\{b^G_{xy}\mid x,y\in W\}$, we may define a \(\Bbbk\)-linear map
\[\Phi:\Sigma_H\longrightarrow e_H\Sigma_Ge_H,\qquad\Phi(b^H_{xy})=b^G_{xy}.\]
This is clearly a $\Bbbk$-isomorphism, and so  \(\Phi\) is an algebra isomorphism.
\end{proof}

\subsection{Projecting decompositions and Mayer--Vietoris sequences}
\label{subsec:mayer-vietoris}

We recall the projecting condition of Hepworth and Willerton
\cite[Definitions~6.2 and 6.3]{HW17} and explain its meaning for distance algebras.  For simplicity, all graphs appearing below are  connected.

\begin{dfn}\label{def:projecting-subgraph}
Let $U$ be a non-empty convex subgraph of a graph $X$.  
For $x\in V(X)$ and $\pi(x)\in V(U)$, we say that $\pi(x)$ is a \defn{projection of $x$ onto $U$} if 
\[
d_X(x,u)
=d_X(x,\pi(x))+d_U(\pi(x),u)
\qquad\text{for all $u\in V(U)$.}
\]
Note that $\pi(x)$ is uniquely determined by this property.
If every $x\in V(X)$ has a projection $\pi(x)\in V(U)$ onto $U$, then we say that $X$ \defn{projects onto $U$}.
\end{dfn}

Since $U$ is convex, we have $d_U(u,v)=d_X(u,v)$ for all
$u,v\in V(U)$.  Notice also that $\pi(u)=u$ for every $u\in V(U)$.
The next result gives a distance algebra interpretation of the
projecting condition.

\begin{prop}\label{prop:projecting-distance-algebra}
Let $U$ be a non-empty convex subgraph of $X$. 
For $x\in V(X)$ and $p\in V(U)$, the following conditions are
equivalent.
\begin{enumerate}
    \item The vertex $p$ is the projection of $x$ onto $U$.
    \item Left-multiplying by $b_{xp}$ 
    \[
    e_p\Sigma_X e_U\longrightarrow e_x\Sigma_Xe_U,
    \qquad a\longmapsto b_{xp}a,
    \]
    defines an isomorphism of right $e_U\Sigma_Xe_U$-modules.
    \item Right-multiplying by $b_{px}$ 
    \[
    e_U\Sigma_X e_p\longrightarrow e_U\Sigma_Xe_x,
    \qquad a\longmapsto ab_{px},
    \]
    defines an isomorphism of left $\Sigma_U$-modules.
\end{enumerate}
Consequently, if $X$ projects onto $U$, then the following hold
\begin{enumerate}[(a)]
\item \[
\Sigma_Xe_U
\cong
\bigoplus_{x\in V(X)}e_{\pi(x)}\Sigma_U
\quad\text{and}\quad 
e_U\Sigma_X
\cong
\bigoplus_{x\in V(X)}\Sigma_Ue_{\pi(x)}
\]
as right and left $e_U\Sigma_Xe_U\cong \Sigma_U$-modules, respectively.

\item $e_U\Sigma_X (1-e_U)\Sigma_X e_U=0$.  In particular, there are algebra isomorphisms \[
\Sigma_U \cong e_U \Sigma_X e_U\cong 
\Sigma_X/\Sigma_X(1-e_U)\Sigma_X.\]
\end{enumerate}
\end{prop}
\begin{proof}
The equivalence between (1) and (2) is straightforward by examining the basis elements; likewise for the equivalence between (1) and (3).  The claim (a) follows immediately from the equivalence.

For (b), consider $u,v\in V(U)$ and $z\in V(X)$.
Recall that $b_{uz}b_{zv}\neq 0$ can occur if and only if $d_X(u,v)=d_X(u,z)+d_X(z,v)$.
Applying projection identities for $u$ and $v$ yields
\begin{align*}
d_X(u,v) &= d_X(u,z)+d_X(z,v) = 2d_X(z,\pi(z))+d_U(\pi(z),u)+d_U(\pi(z),v) \\
& \geq 2d_X(z,\pi(z))+d_U(u,v) \\
& = 2d_X(z,\pi(z))+d_X(u,v),
\end{align*} where the second line follows from the triangle inequality on $U$, and the last follows from the convexity of $U$.
Thus, $b_{uz}b_{zv}=0$ for all $z\notin V(U)$, meaning that $e_U\Sigma_X(1-e_U)\Sigma_Xe_U=0$.
\end{proof}

\begin{rmk}\label{rmk:quotient-not-projecting}
The condition (b) alone does not characterize the projecting condition.  For example, if $U$ is an edge of the cycle $C_3$, then $\Sigma_U\cong\Sigma_{C_3}/\Sigma_{C_3}(1-e_U)\Sigma_{C_3}$, but
$C_3$ does not project onto $U$; see also the discussion following 
\cite[Definition~6.2]{HW17}.
\end{rmk}

\begin{dfn}\label{def:projecting-decomposition}
A \defn{projecting decomposition} is a triple $(X;G,H)$ consisting
of a graph $X$ and subgraphs $G,H\subseteq X$ such that
\begin{enumerate}
    \item $X=G\cup H$;
    \item $G\cap H$ is convex in $X$;
    \item $H$ projects onto $G\cap H$.
\end{enumerate}
\end{dfn}

We next record how such a decomposition is seen inside the distance
algebra of $X$.

\begin{prop}\label{prop:projecting-decomposition-algebra}
Let $(X;G,H)$ be a projecting decomposition and put $K=G\cap H$.
Then $G$ and $H$ are convex in $X$, and there are natural isomorphisms
\[
e_G\Sigma_Xe_G\cong\Sigma_G,
\qquad
e_H\Sigma_Xe_H\cong\Sigma_H,
\qquad
e_K\Sigma_Xe_K
\cong
e_K\Sigma_Ge_K
\cong
e_K\Sigma_He_K
\cong
\Sigma_K.
\]
Moreover, for $h\in V(H)$, $g\in V(G)$, and $p=\pi(h)$, one has
$d_X(h,g)=d_H(h,p)+d_G(p,g).$
\end{prop}

\begin{proof}
There is no edge of $X$ joining a vertex of $G\setminus K$ to a
vertex of $H\setminus K$. Any path passing through $H$ whose
endpoints lie in $G$ begins and ends in $K$. Since $K$ is
convex in $X$, such a path can be replaced by a path in $K$ of
no greater length. Thus $G$ is convex in $X$. Similarly, $H$ is convex in $X$. The corner identifications now follow from
Lemma~\ref{lem:corner-convex-subgraph}. Note that every path from $h$ to $g$ meets $K$. Hence
$d_X(h,g)= \min_{k\in V(K)}
\bigl\{d_H(h,k)+d_G(k,g)\bigr\}$.
For every $k\in V(K)$, the projecting condition and the triangle
inequality in $G$ give
\[
\begin{aligned}
d_H(h,k)+d_G(k,g)
&=
d_H(h,p)+d_K(p,k)+d_G(k,g)\\
&\geq
d_H(h,p)+d_G(p,g).
\end{aligned}
\]
On the other hand, concatenating a shortest path from $h$ to $p$ in
$H$ with a shortest path from $p$ to $g$ in $G$ gives
$d_X(h,g)\leq d_H(h,p)+d_G(p,g)$. The desired distance identity follows.
\end{proof}

We now formulate the Mayer--Vietoris theorem in representation-theoretic terms.

\begin{thm}\textnormal{(cf. \cite[Theorem~6.6]{HW17})}\label{thm:distance-algebra-mv}
Let $(X;G,H)$ be a projecting decomposition and put $K=G\cap H$.
For every $n,\ell\geq0$, the natural inclusions of subgraphs induce split short exact sequences
\[\begin{tikzcd}[row sep=0]
	0 \ar[r] & \Tor_{n,\ell}^K(S,S) \ar[r] & \Tor_{n,\ell}^G(S,S)\oplus \Tor_{n,\ell}^H(S,S)\ar[r]  & \Tor_{n,\ell}^X(S,S)\ar[r] & 0, \\
    0\ar[r] & \Ext^{n,\ell}_X(S,S)\ar[r] & \Ext^{n,\ell}_G(S,S)\oplus\Ext^{n,\ell}_H(S,S)\ar[r] & \Ext^{n,\ell}_K(S,S)\ar[r] & 0.
\end{tikzcd}\]
\end{thm}

Combining this characterization with the Mayer--Vietoris sequence
in Theorem~\ref{thm:distance-algebra-mv} gives the following algebraic form of the gluing property for
diagonal graphs.

\begin{cor}\label{cor:mv-koszul}
Let $(X;G,H)$ be a projecting decomposition and put $K=G\cap H$.
Then the following conditions are equivalent.
\begin{enumerate}
    \item Both $\Sigma_G$ and $\Sigma_H$ are Koszul.
    \item Both $\Sigma_X$ and $\Sigma_K$ are Koszul.
\end{enumerate}
\end{cor}

\subsection{Joins}\label{subsec:join}
\begin{dfn}
    Let $G$ and $H$ be two graphs. The \defn{join} of $G$ and $H$, denoted by $G\vee H$, is the graph with vertex set
$$V(G\vee H)=V(G)\sqcup V(H)$$
and edge set
$$E(G\vee H)=E(G)\sqcup E(H)\sqcup\bigl\{\{x,y\}\mid x\in V(G),\ y\in V(H)\bigr\}.$$
Thus every vertex of $G$ is adjacent to every vertex of $H$.
\end{dfn}

\begin{lem}\label{lem:join-distance}
For $x,x'\in V(G)$,
\[
d_{G\vee H}(x,x')=
\begin{cases}
0, & x=x',\\
1, & \{x,x'\}\in E(G),\\
2, & x\neq x'\text{ and }\{x,x'\}\notin E(G).
\end{cases}
\]
Moreover, $d_{G\vee H}(x,y)=1$ for all $x\in V(G)$ and $y\in V(H)$.
\end{lem}

\begin{proof}
The only point to note is that if two distinct vertices lie in the same
part and are not adjacent there, then they are connected by a path of
length two through any vertex of the other part. The remaining assertions are immediate from the definition of the join.
\end{proof}

\begin{prop}\label{prop:join-matrix}
Let $G$ and $H$ be graphs and consider $\Sigma_{G\vee H}$.
Then via the natural basis-preserving map, we have
$e_G\Sigma_{G\vee H}e_G\cong \Sigma_G$ if and only if $\operatorname{diam}(G)\leq 2$. Therefore, if $\operatorname{diam}(G)\leq 2$ and $\operatorname{diam}(H)\leq 2$, then
\[\Sigma_{G\vee H}\cong\begin{pmatrix}
\Sigma_G & {}_GS_H\\
{}_HS_G &\Sigma_H
\end{pmatrix},\]
where ${}_XS_Y:= S_X\otimes_{\Bbbk} S_Y$ denotes the $\Sigma_X$-$\Sigma_Y$-bimodule given by tensoring their semisimple tops (with left action given by the duality $\omega$ on $X$) for any two graphs $X,Y$.
\end{prop}

\begin{proof}
By Lemma~\ref{lem:join-distance}, the distances in $G\vee H$ agree with the distances in $G$ for all pairs $x,x'\in V(G)$ if and only if $\operatorname{diam}(G)\leq 2$.
Therefore by Lemma \ref{lem:corner-convex-subgraph}, $\Sigma_G\cong e_G\Sigma_{G\vee H}e_G$ if and only if $\operatorname{diam}(G)\leq 2$.
The description of the non-diagonal entries follows from the fact that $e_x \Sigma_{G\vee H} e_y$ is always of rank $1$ over $\Bbbk$.
\end{proof}

\section{Projective dimensions of path-generated right ideals}\label{sec:4-proj-dim}

We keep our assumption of $\Bbbk$ being a unital commutative ring.

\begin{prop}\label{prop:cartan gldim}
The following hold for a finite connected graph $G$.
\begin{enumerate}
    \item The Cartan matrix of $\Sigma_G$ is the all-one matrix.  In particular, it has determinant $0$ if $|V|\ge 2$.
    
    \item The global dimension of $\Sigma_G$ is infinite unless $|V|=1$.
\end{enumerate}
\end{prop}
\begin{proof}
 (1)  Clear from construction.

 (2) We present two proofs.  First note that for an idempotent $e$ of an algebra $A$, the Schur functor $(-)e=-\otimes_AAe\cong \Hom_A(eA,-)$ is exact, and so if a simple $A$-module $S$ has $\operatorname{pdim}(Se)_{eAe}=\infty$, then $\operatorname{pdim}S_{A}=\infty$.  For an edge $\{x,y\}\in E$, we have $\Sigma_{\{x,y\}} = e\Sigma_Ge$ with $e=e_x+e_y$, which is a radical-square-zero Nakayama algebra of infinite global dimension.  Hence, the simple $\Sigma_G$-module $S_x$ associated to $x$ satisfies $\operatorname{pdim}(S_xe)_{\Sigma_{\{x,y\}}}=\infty$, and we are done.

 For a second proof, notice that (1) says that the Cartan matrix is non-invertible (over $\mathbb{Z}$) unless $|V|=1$, and so the claim follows from the fact that an algebra of finite global dimension must have an invertible Cartan matrix. The proof of this fact is a linear algebra exercise by considering  modules as elements of the Grothendieck group and their projective resolutions; we give a demonstration in our special case, which has an even easier computation, for completeness.

We note the following lemma.
\begin{lem}\label{lem:fin pdim module dimvec}
Let $M$ be a $\Sigma_G$-module with finite projective dimension. 
Then the dimension vector of $M$ has constant entries, i.e., $\sdim M\in \mathbb{Z}(1,1,\ldots,1)$.
\end{lem}
\begin{proof}
The dimension vector of $M$ satisfies $\sdim M = \sum_{r=0}^m(-1)^r \sdim P_r$, where
$$0\to P_m\to \cdots \to P_1\to P_0\to M\to 0$$
is a projective resolution of $M$.  In other words, $\sdim M \in \mathbb{Z}\text{-span} \{\sdim P_x\}_{x\in V} = \mathbb{Z}(1,1,\ldots, 1)$.    
\end{proof}

Suppose that $|V|\geq 2$.
Let $S_x$ be the simple module associated to $x\in V$.  Then its dimension vector is the $x$-th standard basis vector $\sdim S_x=(0,\dots,0,1,0,\dots,0)\notin \mathbb Z(1,1,\dots,1)$.  Hence, $S_x$ has infinite projective dimension, and so $\Sigma_G$ has infinite global dimension.
\end{proof}
\begin{rmk}\label{rmk:mayer-vietoris}
    Part (2) can also be shown using the Mayer--Vietoris sequence of magnitude (co)homology \cite{BK21}, which is essentially the same argument as the first proof of this claim.
\end{rmk}

Consider the cyclic right $\Sigma_G$-module (right principal ideal) $P_x^y:=b_{xy}\Sigma_G$ generated by $b_{xy}$. 
Clearly, we have $P_x^y\subseteq P_x^x=:P_x$ for all $x,y\in V$.
By the multiplication rule of $\Sigma_G$, the module $P_x^y$ has a $\Bbbk$-basis
$$\{b_{xz}\mid z\in V,\ d(x,z)=d(x,y)+d(y,z)\}.$$
Equivalently, $b_{xz}$ belongs to $P_x^y$ precisely when $y$ lies on a shortest path from $x$ to $z$.

\begin{lem}\label{lem:syzygy}
Let $x,y\in V$ be adjacent vertices. Then the projective cover of $P_x^y$ is given by
$$e_y\Sigma_G=P_y\rightarrow P_x^y,\quad a\mapsto b_{xy}a$$
and the syzygy is given by
$$\Omega_{\Sigma_G}(P_x^y)=\operatorname{span}_{\Bbbk}\{b_{yz}\mid d(x,z)\leq d(y,z)\}.$$
Moreover, $\Omega_{\Sigma_G}(P_x^y)$ is always isomorphic to a submodule of $P_y^x$, with equality when  there is no vertex $z\in V$ such that $d(x,z)=d(y,z)$.  This condition holds, for example, when $G$ is bipartite. Hence, in particular, it holds when $G$ is a tree.
\end{lem}

\begin{proof}
Since $x$ and $y$ are adjacent, we have $d(x,y)=1$. The module $P_x^y=b_{xy}\Sigma_G$ is generated by $b_{xy}$. Hence we have a natural surjective homomorphism
$\pi:e_y\Sigma_G\rightarrow P_x^y$ given by $e_y\mapsto b_{xy}$.
Equivalently, $\pi(b_{yz})=b_{xy}b_{yz}$ for all $z\in V$. 
Since $x$ and $y$ are adjacent, the triangle inequality gives $d(x,z)\leq 1+d(y,z)$ for every $z\in V$.
By the multiplication rule in $\Sigma_G$, we get that
\begin{align*}
\Omega_{\Sigma_G}(P_x^y)\cong\mathrm{Ker}\pi&=\operatorname{span}_{\Bbbk}\{b_{yz}\mid d(x,z)\neq 1+d(y,z)\}\\
&=\operatorname{span}_{\Bbbk}\{b_{yz}\mid d(x,z)\leq d(y,z)\}. 
\end{align*}

On the other hand, since $P_y^x$ has a $\Bbbk$-basis $\{b_{yz}\mid d(y,z)=1+d(x,z)\}$, we have  $\Omega_{\Sigma_G}(P_x^y)\cong P_y^x$ precisely when there is no vertex $z$ with $d(x,z)=d(y,z)$.
\end{proof}

\begin{prop}\label{prop:pdim Pxy}\label{prop:gldim-infity}
    $P_x^y$ has infinite projective dimension if $x\neq y$.
\end{prop}

\begin{proof}
Assume $x\neq y$. The dimension vector is $\sdim P_x^y=(\epsilon_z)_{z\in V}$, where
$$
\epsilon_z=
\begin{cases}
1, & \text{if } d(x,z)=d(x,y)+d(y,z),\\
0, & \text{otherwise}.
\end{cases}
$$
In particular, $\epsilon_y=1$, while $\epsilon_x=0$, since $x\neq y$. Hence $\sdim P_x^y$ is not a multiple of $(1,1,\dots,1)$. Therefore by Lemma \ref{lem:fin pdim module dimvec}, $\operatorname{pdim}_{\Sigma_G} P_x^y=\infty$ for $x\neq y$.
\end{proof}

Recall that \defn{Han's conjecture} \cite{Han06} predicts that, for a finite-dimensional algebra $A$ over a field $\Bbbk$, the Hochschild homology $\mathrm{HH}_n(A):=\operatorname{Tor}_n^{A^{\operatorname{op}}\otimes_\Bbbk A}(A,A)=0$ for all $n\gg 0$ if and only if the global dimension of $A$ is finite.  Note that the if-direction always holds \cite[Proposition 3.3]{Han06}.

\begin{prop}
Suppose that $\Bbbk$ is a field and $|V|\ge 2$.  Then $\mathrm{HH}_n(\Sigma_G)\neq 0$ for all $n\geq 0$.  In particular, $\Sigma_G$ satisfies Han's Hochschild homology conjecture.
\end{prop}
\begin{proof}
Any edge $\{x,y\}\in E$ defines a subalgebra $\Sigma_{\{x,y\}} = e\Sigma_Ge$ with $e=e_x+e_y$, which is a radical-square-zero Nakayama algebra of infinite global dimension.
Moreover, this subalgebra can be extended to a unital split subalgebra $B:=e\Sigma_Ge\oplus \Bbbk(1-e)\subseteq \Sigma_G$, where $B$ is an algebra retract of $\Sigma_G$.  
Now the non-vanishing of Hochschild homology groups follows from the functoriality of Hochschild homology \cite{Lod98,BM10}. 
\end{proof}
\begin{rmk}
While we believe the analogous result for cohomology should hold, we do not have any solid supporting evidence yet.
\end{rmk}


\section{Biseriality and representation-finiteness}\label{sec:biserial}

From now on until the end of the article, we take $\Bbbk$ to be a field.

\begin{dfn}
Let $A$ be a finite-dimensional algebra.
\begin{enumerate}
\item The algebra $A$ is called \defn{multiserial} if, for every indecomposable projective right or left $A$-module $P$, there exist, for some $m\geq0$, uniserial submodules $U_1,\ldots,U_m\subseteq\operatorname{rad}P$ such that
$$
\operatorname{rad}P=U_1+\cdots+U_m
$$
and $U_i\cap U_j$ is either zero or simple whenever $i\neq j$.

\item The algebra $A$ is called \defn{biserial} if one can always take
$m=2$ in $(1)$, allowing one or both of $U_1,U_2$ to be zero.

\item Let $A=\Bbbk Q/I$ be a bound quiver algebra. It is called \defn{special multiserial} if it satisfies the following condition:
\begin{itemize}
\item[(SM)] for every arrow $\alpha\in Q_1$, there is at most one arrow $\beta\in Q_1$ such that $\alpha\beta\notin I$, and at most one arrow $\gamma\in Q_1$ such that $\gamma\alpha\notin I$.
\end{itemize}
It is called \defn{special biserial} if, in addition, it satisfies:
\begin{itemize}
\item[(SB)] for every vertex $i\in Q_0$, there are at most two arrows starting at $i$ and at most two arrows ending at $i$.
\end{itemize}

\item The algebra $A=\Bbbk Q/I$ is called a \defn{string algebra} if it is special biserial and $I$ is generated by paths.

\item The algebra $A=\Bbbk Q/I$ is called a \defn{gentle algebra} if it is special biserial and satisfies the following additional conditions:
\begin{itemize}
\item[(G1)] for every arrow $\alpha\in Q_1$, there is at most one arrow $\beta\in Q_1$ such that $\alpha\beta\in I$, and at most one arrow $\gamma\in Q_1$ such that $\gamma\alpha\in I$;
\item[(G2)] the ideal $I$ is generated by paths of length two.
\end{itemize}
\end{enumerate}
\end{dfn}

We have the following implications:
$$\begin{tikzcd}[column sep=small,row sep=small]
&&& \text{biserial} \arrow[dr, Rightarrow] & \\
\text{gentle} \arrow[r, Rightarrow]
& \text{string} \arrow[r, Rightarrow]
& \text{special biserial}
  \arrow[ur, Rightarrow]
  \arrow[dr, Rightarrow]
&& \text{multiserial} \\
&&& \text{special multiserial} \arrow[ur, Rightarrow] &
\end{tikzcd}$$
For these definitions and further related results, we refer the reader to \cite{Fuller68,SW83,GS16}.

\begin{eg}\label{eg:line and cycle}
\begin{enumerate}
\item Consider the path graph $G=A_n$. Then $\Sigma_G$ is given by the quiver
$$
\begin{tikzcd}
	1 & 2 & \cdots & n
	\arrow["{a_1}", bend left, from=1-1, to=1-2]
	\arrow["{b_1}", bend left, from=1-2, to=1-1]
	\arrow["{a_2}"{pos=0.6}, bend left, from=1-2, to=1-3]
	\arrow["{b_2}"{pos=0.4}, bend left, from=1-3, to=1-2]
	\arrow["{a_{n-1}}"{pos=0.45}, bend left, from=1-3, to=1-4]
	\arrow["{b_{n-1}}"{pos=0.55}, bend left, from=1-4, to=1-3]
\end{tikzcd}
$$
with relations $a_ib_i$ and $b_ia_i$ for all $1\leq i<n$. Thus $\Sigma_{A_n}$ is gentle.

\item Consider the cycle $G=C_n$ with $n\geq3$. The double quiver
$Q_G$ satisfies (SB), and every arrow has at most one successor and at
most one predecessor whose product is non-zero. Hence $\Sigma_{C_n}$ is
special biserial. More precisely, $\Sigma_{C_n}$ is a string algebra
when $n$ is odd, whereas it is a self-injective special biserial algebra when $n$ is even.

\item Let $G=K_n$ be a complete graph. No path of length two in $Q_G$ is a shortest path between its endpoints. Consequently, every product of two arrows vanishes and
$$\Sigma_{K_n}\simeq\Bbbk Q_{K_n}/J^2,
\qquad\operatorname{rad}^2\Sigma_{K_n}=0,$$
where $J$ denotes the arrow ideal of $\Bbbk Q_{K_n}$. In particular, condition (SM) holds, so $\Sigma_{K_n}$ is special multiserial for every $n\geq1$.
\end{enumerate}
\end{eg}

For a vertex $x\in V(G)$, let $N(x)$ denote its \defn{open neighborhood} and let $N[x]:=N(x)\cup\{x\}$ denote its \defn{closed neighborhood}. Set $\overline N(x):=V(G)\setminus N[x]$, the set of \defn{non-neighbors} of $x$ in $G$.

\begin{prop}\label{prop:multiserial}
Let $G=(V,E)$ be a connected graph. The following are equivalent:
\begin{enumerate}
\item $\Sigma_G$ is multiserial;
\item $\Sigma_G$ is special multiserial;
\item for every edge $\{x,y\}\in E$, one has
$\left|N(x)\setminus N[y]\right|\leq1$, and $\left|N(y)\setminus N[x]\right|\leq1$;
\item $G$ is a path graph $A_n$, a cycle $C_n$, or there are pairwise distinct vertices $x_1,y_1,\ldots,x_r,y_r\in V$ for some $r\geq0$ such that
$$E=\binom{V}{2}\setminus
\bigl\{\{x_1,y_1\},\ldots,\{x_r,y_r\}\bigr\}.$$
\end{enumerate}
\end{prop}

\begin{proof}
The implication $(2)\Rightarrow(1)$ is immediate.

We prove $(3)\Rightarrow(2)$. Let $a_{xy}:x\to y$ be an arrow of $Q_G$, and let $z\in N(y)$. By the multiplication rule in the distance algebra,
$$a_{xy}a_{yz}\notin I_G
\quad\Longleftrightarrow\quad
d(x,z)=2
\quad\Longleftrightarrow\quad
z\in N(y)\setminus N[x].$$
Likewise, for an arrow $a_{wx}:w\to x$, one has
$a_{wx}a_{xy}\notin I_G$ if and only if $w\in N(x)\setminus N[y]$.
Thus the two inequalities in $(3)$ say precisely that every arrow has at most one successor and at most one predecessor whose product is non-zero. Hence $\Sigma_G$ satisfies (SM).

For $(1)\Rightarrow(3)$, suppose that there are an edge $\{x,y\}$ and
distinct vertices $z_1,z_2\in N(y)\setminus N[x]$. Then
$b_{xy}b_{yz_i}=b_{xz_i}\neq0$ for $i=1,2$. No uniserial submodule of $\operatorname{rad}P_x$ can contain $b_{xy}$. Indeed, such a submodule would contain both $b_{xz_1}$ and $b_{xz_2}$, and hence the two non-zero submodules $b_{xz_1}\Sigma_G$ and $b_{xz_2}\Sigma_G$. Note that these submodules are incomparable. On the other hand, if $\operatorname{rad}P_x$ were a sum of uniserial submodules $U_1,\ldots,U_m$, then $(\operatorname{rad}P_x)e_y=\sum_{j=1}^m U_je_y=\Bbbk b_{xy}$. Since $U_je_y\subseteq P_xe_y=\Bbbk b_{xy}$, this equality shows that one of the uniserial submodules $U_j$ contains $b_{xy}$, a contradiction. Therefore $\Sigma_G$ is not multiserial. Interchanging $x$ and $y$ proves both inequalities in $(3)$.

It remains to prove $(3)\Longleftrightarrow(4)$. Suppose first that $(3)$ holds. If every vertex has valency at most two, then connectedness implies that $G$ is a path graph or a cycle. Now suppose that $\val(v)\geq3$ for some vertex $v$, and put
$$X:=N(v),\qquad Y:=V\setminus N[v].$$
For every $u\in X$, condition $(3)$ applied to $\{u,v\}$ shows that $u$ is non-adjacent to at most one vertex of $X$ and is adjacent to at most one vertex of $Y$.

Let $w\in Y$ satisfy $d(v,w)=2$, and choose $u\in X\cap N(w)$. Since $v\in N(u)\setminus N[w]$,
condition $(3)$ applied to $\{u,w\}$ shows that every neighbor of $u$ other than $v$ and $w$ is adjacent to $w$. We claim that $w$ is adjacent to every vertex of $X$. Otherwise, choose $s\in X$ not adjacent to $w$. Then $s$ is not adjacent to $u$. Choose $t\in X\setminus\{u,s\}$. Since every vertex of $X$ has at most one non-neighbor in $X$, the vertex $t$ is adjacent to both $u$ and $s$; hence it is also adjacent to $w$. But then $v,s\in N(t)\setminus N[w]$, contrary to $(3)$ applied to $\{t,w\}$. Thus $X\subseteq N(w)$.

There is no vertex at distance three from $v$. Indeed, if $v,u,w,q$ were a shortest path, then $X\subseteq N(w)\setminus N[q]$,
contrary to $|X|=\val(v)\geq3$ and condition $(3)$ applied to $\{w,q\}$. By connectedness, every vertex of $Y$ therefore has distance two from $v$ and is adjacent to every vertex of $X$. For any $u\in X$, it follows that $Y=N(u)\setminus N[v]$,
so $(3)$ gives $|Y|\leq1$. Consequently, every vertex of $G$ has at most one non-neighbor. Thus the non-edges of $G$ are pairwise vertex-disjoint.

Conversely, path graphs and cycles clearly satisfy $(3)$. In the last case of $(4)$, every vertex has at most one non-neighbor, so both inequalities in $(3)$ hold. This completes the proof.
\end{proof}

As a special case, we recover the biserial distance algebras.

\begin{prop}\label{prop:biserial}
Let $G=(V,E)$ be a connected graph. The following are equivalent:
\begin{enumerate}
\item $\Sigma_G$ is biserial;
\item $\Sigma_G$ is special biserial;
\item $G$ is either the path graph $A_n$ with $n\geq1$ or a cycle $C_n$ with $n\geq3$.
\end{enumerate}
\end{prop}

\begin{proof}
The implication $(3)\Rightarrow(2)$ follows from Example~\ref{eg:line and cycle}, and $(2)\Rightarrow(1)$ is immediate.

For $(1)\Rightarrow(3)$, the number of arrows starting at a vertex $x$ of $Q_G$ is $\val(x)$, and
$$\operatorname{rad}P_x/\operatorname{rad}^2P_x
\simeq\bigoplus_{y\in N_G(x)}S_y.$$
Since $\operatorname{rad}P_x$ is a sum of at most two uniserial modules, its top has composition length at most two. Therefore $\val(x)\leq2$ for every $x\in V$. Connectedness now implies that $G$ is a path graph or a cycle.
\end{proof}

\begin{rmk}
   Note that a distance algebra $\Sigma_G$ is quasi-biserial \cite{Xin24} if and only if it is biserial. Indeed, suppose that $\Sigma_G$ is not biserial. There is a vertex $y$ of valency at least three. At this vertex $y$, there are at least three incoming arrows and at least three outgoing arrows.
\end{rmk}

\begin{dfn}
Let $A$ be a finite-dimensional $\Bbbk$-algebra.
We say that $A$ is of \defn{tame} representation type, or simply \emph{tame}, if, for every positive integer $d\geq 1$, there exist finitely many $A$--$\Bbbk[x]$-bimodules $M_1,\ldots,M_s$, free of finite rank as left $\Bbbk[x]$-modules, such that all but finitely many isomorphism classes of indecomposable $A$-modules of dimension $d$ are of the form $S\otimes_{\Bbbk[x]}M_i$, where $S$ is a simple $\Bbbk[x]$-module. Denote by $\mu_A(d)$ the least possible such $s$ -- the number of one-parameter families of indecomposable modules of dimension $d$.
We say that
\begin{itemize}
\item $A$ is \defn{domestic} if it is tame and there is some non-negative integer $C$ such that $\mu_A(d)\le C$ for all $d$; and

\item $A$ is \defn{of polynomial growth} if it is tame and there is some integer $N\ge 1$ such that $\mu_A(d)\leq d^N$ for all $d\geq 1$.
\end{itemize}
\end{dfn}

\begin{rmk}
Representation-finite algebras are also tame; we abbreviate the case when it is representation-infinite and tame by \defn{infinite-tame}.
By convention in the literature, when we say that an algebra is \emph{not} of polynomial growth, then it is ``beyond'' polynomial growth, i.e., infinite-tame and non-domestic.  
\end{rmk} 

Let us recall some standard facts. 
Every biserial algebra is tame \cite{CB95}. 
The indecomposable non-projective modules of a special biserial algebra belong to one of the two families called string and band modules.
The one-parameter families occur only in band modules, and each family is described by a combinatorial gadget called \defn{band}; we refer the reader to \cite{BR87} for details.
In particular, a special biserial algebra is 
\begin{itemize}
    \item representation-finite if and only if  there is no band;
    \item infinite-tame domestic if and only if there are only finitely many bands.
\end{itemize}
Moreover, domesticity is equivalent to polynomial growth for string algebras; see, for example, in \cite{BR87} or \cite[Corollary~B]{EKM26}, and, for self-injective special biserial algebras, see \cite[Theorems~2.1 and~2.2]{ES92}.

\begin{prop}\label{prop:biserial rep type}
The following hold.
\begin{enumerate}
    \item $\Sigma_{A_n}$ is representation-finite for all $n\ge 1$.
    \item $\Sigma_{C_3}$ and $\Sigma_{C_4}$ are representation-domestic.
    \item $\Sigma_{C_n}$ is representation-tame and not of polynomial growth for all $n\geq 5$.
\end{enumerate}
\end{prop}
\begin{proof}
(1) Let $G=A_n$. Then $\Sigma_G$ is a gentle algebra (Example \ref{eg:line and cycle} (1)).  Suppose on the contrary that $\Sigma_G$ is representation-infinite.  Then there is a band, which means that there is some $1\le i\le n$ and a subword (up to inversion) of the form $ab^{-1}$ with $t(a)=t(b)$.  But this is only possible when $b=a$, violating the definition of a string.

(2)  
For $C_3$, up to cyclic permutation and inversion, its unique band is represented by
$$a_{x_1x_2}a_{x_3x_2}^{-1}a_{x_3x_1}a_{x_2x_1}^{-1}a_{x_2x_3}a_{x_1x_3}^{-1}.$$

For $C_4$, up to cyclic permutation and inversion, the bands are represented by
$$a_{x_1x_2}a_{x_3x_2}^{-1}a_{x_3x_4}a_{x_1x_4}^{-1},\qquad a_{x_2x_1}^{-1}a_{x_2x_3}a_{x_4x_3}^{-1}a_{x_4x_1}.$$
Indeed, there are exactly two such bands. Therefore $\Sigma_{C_3}$ and $\Sigma_{C_4}$ are domestic.

(3)  Label the vertices of $C_n$ by $x_i$ for $i\in\mathbb Z/n\mathbb Z$, and write
$$\alpha_i=a_{x_ix_{i+1}}:x_i\to x_{i+1},\qquad\beta_i=a_{x_{i+1}x_i}:x_{i+1}\to x_i$$
for the two arrows corresponding to the edge $\{x_i,x_{i+1}\}$, where indices are taken modulo $n$.
We construct bands by moving clockwise around the cycle. At the $i$-th step, there are two possible letters: $\alpha_i$ or $\beta_i^{-1}$.
Thus $\alpha_i$ means that the next composition is direct, while $\beta_i^{-1}$ means that the next composition is inverse. In both cases the word moves from $x_i$ to $x_{i+1}$. For example, the following two cyclic words are bands:
$$P=\alpha_0\beta_1^{-1}\alpha_2\beta_3^{-1}\cdots \alpha_{2n-2}\beta_{2n-1}^{-1},$$
and
$$Q=\alpha_0\alpha_1\beta_2^{-1}\alpha_3\beta_4^{-1}\cdots
\alpha_{2n-3}\beta_{2n-2}^{-1}\beta_{2n-1}^{-1},$$
where all indices are taken modulo $n$. More generally, for every $m\geq 1$, the cyclic word $P^mQ$ is a band. Indeed, $P^mQ$ are pairwise distinct up to cyclic permutation and inversion, since they have different lengths. They are also primitive, because each cyclic word $P^mQ$ contains exactly one occurrence of two consecutive direct compositions. Hence $\Sigma_{C_n}$ has infinitely many primitive bands, and therefore $\Sigma_{C_n}$ is not domestic, and not of polynomial growth.
\end{proof}



\begin{thm}\label{thm:rep-finite}
$\Sigma_G$ is representation-finite if and only if $G$ is a path graph $A_n$.
\end{thm}

\begin{proof}
The if-direction is shown in Proposition \ref{prop:biserial rep type}.
Suppose now that $G$ is not a path graph.

Recall first that, for every ideal $I$ of $A$, the category $\operatorname{mod} A/I$ identifies with the full subcategory of $\operatorname{mod} A$ consisting of modules annihilated by $I$. In particular, representation-finiteness is inherited by quotient algebras. It is also inherited by idempotent truncations: if $A$ is
representation-finite and $e\in A$ is an idempotent, then $eAe$ is
representation-finite. Indeed, the induction functor
$$-\otimes_{eAe}eA:\operatorname{mod}eAe\longrightarrow
\operatorname{mod}A$$
is fully faithful. Thus $\operatorname{mod}eAe$ identifies with a full subcategory of $\operatorname{mod}A$. Thus, it suffices to find a representation-infinite quotient of an idempotent truncation of $\Sigma_G$.

Suppose first that $G$ contains a cycle. Choose a cycle $C$ of minimal
length in $G$, and set $e_C=\sum_{x\in V(C)}e_x$. The minimality of $C$ implies that $C$ is a convex full subgraph of $G$. Hence, by
Lemma~\ref{lem:corner-convex-subgraph}, $e_C\Sigma_Ge_C\simeq \Sigma_C$. Now consider the quotient
$$B_C:=\Sigma_C/\operatorname{rad}^2\Sigma_C\simeq \Bbbk Q_C/J_C^2,$$
where $J_C$ denotes the arrow ideal of $\Bbbk Q_C$.
By stable equivalence, $\Bbbk Q/J^2$ is representation-finite if and
only if the separated quiver $Q^{\operatorname{sep}}$ is a disjoint
union of Dynkin quivers; see \cite[Section~X.2]{ARS95}. Since $C$ is a
cycle, the underlying graph of $Q_C^{\operatorname{sep}}$ contains a
cycle. Consequently, $Q_C^{\operatorname{sep}}$ is not a disjoint union of Dynkin quivers, and hence $B_C$ is representation-infinite. It follows that if $\Sigma_G$ is representation-finite, then $G$ must be a tree.

Now assume that $G=T$ is a tree. If $T$ is not a path, then $T$ has a vertex of valency at least $3$. Choose three distinct neighbors of such a vertex. The full subtree on these four vertices is the star $D_4$. Let $e$ be the sum of the primitive idempotents corresponding to these four vertices. Then by Lemma \ref{lem:corner-convex-subgraph}, $e\Sigma_Te\simeq \Sigma_{D_4}$. Thus it is enough to show that $\Sigma_{D_4}$ is representation-infinite. 

For every integer $n\geq1$, let $J_n$ denote the nilpotent Jordan block of size $n$.  Consider the representation $N_n$ given by
\[
\begin{tikzcd}[column sep=large]
& & \Bbbk^{2n}
\arrow[dd, "\BTwoN", shift left=0.8ex]
& & \\[1em]
& & & & \\[1em]
\Bbbk^{2n}
\arrow[rr, "\BOneN"]
& &
\Bbbk^{3n}
\arrow[ll, "\AOneN", shift left=1.2ex]
\arrow[uu, "\ATwoN", shift left=0.8ex]
\arrow[rr, "\AThreeN"]
& &
\Bbbk^{2n}
\arrow[ll, "\BThreeN", shift left=1.2ex]
\end{tikzcd}
\]
where $I_n$ is the $n\times n$ identity matrix and every $0$ in the
displayed block matrices denotes the $n\times n$ zero matrix.  Denote
the three pairs of displayed maps by $B_i^{(n)}$ and $A_i^{(n)}$,
respectively.  A direct calculation gives
\[
A_i^{(n)}B_i^{(n)}
=
B_i^{(n)}A_i^{(n)}
=0
\qquad (i=1,2,3).
\]
Therefore, $N_n$ is a $\Sigma_{D_4}$-module.

We first show that $N_n$ is indecomposable.  Let $(X,Y_1,Y_2,Y_3)$
be an endomorphism of $N_n$, where $X$ is the map at the central
vertex and the $Y_i$ are the maps at the three outer vertices.  The
commutativity conditions
\[
XB_i^{(n)}=B_i^{(n)}Y_i,
\qquad
Y_iA_i^{(n)}=A_i^{(n)}X
\]
give, by direct matrix calculation,
\[
X=\operatorname{diag}(P,P,P),
\qquad
Y_i=
\begin{pmatrix}
P&0\\
C_i&P
\end{pmatrix}
\quad (i=1,2,3),
\qquad
PJ_n=J_nP,
\]
for some matrices $P$ and $C_i$. 
Since $J_n$ is a single Jordan block, every matrix commuting with
$J_n$ is a polynomial in $J_n$.  Since $J_n^n=0$, we have $
P\in\Bbbk[t]/(t^n)$.
Suppose now that the endomorphism $(X,Y_1,Y_2,Y_3)$ is idempotent.  Since an 
idempotent in $\Bbbk[t]/(t^n)$ is either the identity or zero, we have $P=0$ or $P=I_n$.  If $P=0$, the
endomorphism has square zero, and hence it must be zero.  On the other hand, if $P=I_n$, then $Y_i^2=Y_i$ implies that $C_i=0$, and hence the endomorphism is the identity.  Thus $N_n$ has no non-trivial
idempotent endomorphisms, which means that $N_n$ is indecomposable.

Finally, if $n\neq m$, then $N_n\not\cong N_m$. Therefore, $\Sigma_{D_4}$ has infinitely many pairwise non-isomorphic indecomposable modules over every field, and is representation-infinite.
\end{proof}


\section{Self-injectivity}\label{sec:self-inj}

\subsection{Self-injectivity is interval-antipodal}
In what follows, we characterize self-injectivity of $\Sigma_G$ in terms of $G$.

\begin{dfn}\label{def:antipodal}
Let $G=(V,E)$ be a graph. \begin{enumerate}
    \item For $x,y\in V$, define the \defn{interval} between $x$ and $y$ by
$$I(x,y):=\{z\in V\mid d(x,y)=d(x,z)+d(z,y)\}.$$
Thus $z\in I(x,y)$ if and only if $z$ lies on a shortest path from $x$ to $y$.
\item For $x\in V$, define
$$\Max(x):=\{y\in V\mid d(x,z)\leq d(x,y)\text{ for every vertex }z\text{ adjacent to }y\}.$$
We call $\Max(x)$ the set of \defn{local maxima} of the distance function $d(x,-)$.
\item We say that $G$ is \defn{interval-antipodal} if, for every $x\in V$, there exists a unique vertex $\nu(x)\in V$ such that $I(x,\nu(x))=V$.
In this case, $\nu(x)$ is called the \defn{antipode} of $x$.
\end{enumerate}
\end{dfn}

\begin{eg}\label{ex:interval-antipodals}
(1) Every even cycle $C_{2n}$ is naturally interval-antipodal with antipode $\nu(i) = i+n$ modulo $2n$ for all $1\le i\le 2n$.

(2) The $n$-dimensional hypercube $Q_n$ is interval-antipodal.  Indeed, we have
$$d(X,Y)+d(Y,X^{\mathrm c})
=|X\mathbin{\triangle}Y|
 +|Y\mathbin{\triangle}X^{\mathrm c}|
=n=d(X,X^{\mathrm c}).$$
Thus $I(X,X^{\mathrm c})=V(Q_n)$. Moreover, $X^{\mathrm c}=[n]\setminus X$ is the unique vertex at distance $n$ from $X$, so it is the unique antipode of $X$.
By the same argument, the doubly punctured $n$-dimensional hypercube $Q_n^{\pm}$ is also interval-antipodal with the same antipodal involution.

(3) Recall that the vertices of the Johnson graph $J(2n,n)$ are the $n$-subsets of $[2n]$, with distance $d(X,Y)=n-|X\cap Y|$.  Similarly to (2), $J(2n,n)$ is interval-antipodal with the antipode given by $X\mapsto X^{\mathrm c}=[2n]\setminus X$.  For example, the six-vertex Johnson graph $J(4,2)$ can be drawn as follows:
\[
\begin{tikzpicture}[
    scale=0.70,
    vertex/.style={
        circle,
        draw,
        fill=white,
        minimum size=10mm,
        inner sep=1pt,
        outer sep=0pt,
        font=\small
    }
]
\draw (0,0) circle[radius=3.2];

\node[vertex] (v12) at (90:3.2)  {$\{1,2\}$};
\node[vertex] (v13) at (210:3.2) {$\{1,3\}$};
\node[vertex] (v14) at (330:3.2) {$\{1,4\}$};

\node[vertex] (v34) at (270:1.25) {$\{3,4\}$};
\node[vertex] (v24) at (30:1.25)  {$\{2,4\}$};
\node[vertex] (v23) at (150:1.25) {$\{2,3\}$};

\draw (v23) -- (v24);
\draw (v24) -- (v34);
\draw (v34) -- (v23);

\draw (v12) -- (v23);
\draw (v12) -- (v24);

\draw (v13) -- (v23);
\draw (v13) -- (v34);

\draw (v14) -- (v24);
\draw (v14) -- (v34);
\end{tikzpicture}
\]
\end{eg}

\begin{lem}\label{lem:proj-inj-max}
Let $x\in V$. Then the following conditions are equivalent.
\begin{enumerate}
\item The projective module $P_x=e_x\Sigma_G$ is injective.
\item The set $\Max(x)$ consists of exactly one vertex.
\item There exists a vertex $y\in V$ such that $I(x,y)=V$.
\end{enumerate}
Moreover, if these conditions hold and $\Max(x)=\{y\}$, then
$P_x\cong I_y$, where $I_y$ is the injective envelope of $S_y$.
\end{lem}

\begin{proof}
The module $P_x$ has $\Bbbk$-basis $\{b_{xz}\mid z\in V\}$. Since the radical of $\Sigma_G$ is generated by the arrows, an element of $P_x$ lies in $\operatorname{soc}P_x$ if and only if it is annihilated by every arrow. Hence
\[
\operatorname{soc}P_x=
\bigoplus_{y\in \Max(x)}\Bbbk b_{xy}
\cong
\bigoplus_{y\in \Max(x)}S_y .
\]
Since $P_x$ is finite-dimensional, its socle is essential. Therefore the injective envelope of $P_x$ is
\[
I(P_x)\cong \bigoplus_{y\in \Max(x)} I_y .
\]
Since $\dim_{\Bbbk}P_x=|V|$ and
$\dim_{\Bbbk}I_y=\dim_{\Bbbk}\Sigma_Ge_y=|V|$ for every $y\in V$,
the module $P_x$ is injective if and only if $\Max(x)$ consists of exactly one vertex. In this case, if $\Max(x)=\{y\}$, then $P_x\cong I_y$.

It remains to show the equivalence between (2) and (3). Suppose first that
$\Max(x)=\{y\}$. Let $v\in V$. If $v=y$, there is nothing to prove. If $v\neq y$, then $v\notin \Max(x)$. Hence there exists a vertex $v_1$ adjacent to $v$ such that
$d(x,v_1)=d(x,v)+1$. Repeating this argument, and using that $G$ is finite, we see that the process must stop at a local maximum, hence at $y$. Thus we get a path
\[
v-v_1-\cdots-v_r=y
\]
such that $d(x,v_s)=d(x,v)+s$ for all $s$. Thus
$d(x,y)=d(x,v)+r$. Since $d(v,y)\leq r$, the triangle inequality gives
\[
d(x,y)\leq d(x,v)+d(v,y)\leq d(x,v)+r=d(x,y).
\]
Therefore $d(x,y)=d(x,v)+d(v,y)$, and hence $v\in I(x,y)$. Since $v$ was arbitrary, $I(x,y)=V$.

Conversely, suppose that $I(x,y)=V$. First, $y\in \Max(x)$. Indeed, if
$z$ is adjacent to $y$, then $z\in I(x,y)$, and so
\[
d(x,y)=d(x,z)+d(z,y)=d(x,z)+1.
\]
Thus no neighbor of $y$ has distance $d(x,y)+1$ from $x$. Now let
$v\neq y$. Since $v\in I(x,y)$, there is a shortest path from $x$ to $y$ passing through $v$. Let $z$ be the next vertex after $v$ on this path toward $y$. Then $d(x,z)=d(x,v)+1$, and hence $v\notin \Max(x)$. Therefore $\Max(x)=\{y\}$.
\end{proof}

\begin{cor}\label{cor:proj-inj-exists}
The algebra $\Sigma_G$ has an indecomposable projective-injective module if and only if there exist vertices $x,y\in V$ such that $I(x,y)=V$. In this case, we have $P_x=I_y$ and $P_y=I_x$.
\end{cor}

\begin{proof}
This follows immediately from Lemma~\ref{lem:proj-inj-max}.
\end{proof}

\begin{eg}
Let $G=A_n$ be the path graph with $n$ vertices.  Then we have $I(1,n)=V(G)$.  Hence, $P_1\cong I_n$ and $P_n\cong I_1$.
\end{eg}

\begin{lem}\label{lem:interval-antipodal-involution}
For a graph $G$, the following conditions are equivalent.
\begin{enumerate}
\item $G$ is interval-antipodal.
\item There exists a permutation $\nu:V\to V$ such that
$I(x,\nu(x))=V$ for every $x\in V$.
\end{enumerate}
In this case, the permutation $\nu$ is involutive, uniquely determined, and has no fixed point unless $|V|=1$.
\end{lem}

\begin{proof}
Assume first that $G$ is interval-antipodal. Then, for every $x\in V$, there exists a unique vertex $\nu(x)\in V$ such that $I(x,\nu(x))=V$.
This defines a map $\nu:V\to V$. Since
$$I(x,\nu(x))=I(\nu(x),x)=V,$$
the uniqueness of the antipode of $\nu(x)$ implies that $\nu^2(x)=x$ for every $x\in V$. Hence $\nu^2=\operatorname{id}_V$, and in particular $\nu$ is a permutation.
Conversely, suppose that there exists a permutation $\nu:V\to V$ such that $I(x,\nu(x))=V$ for every $x\in V$. Then for every $x\in V$ there exists a vertex $y=\nu(x)$ such that $I(x,y)=V$. By Lemma~\ref{lem:proj-inj-max}, such a vertex is unique. Hence $G$ is interval-antipodal. 
The same uniqueness also shows that the permutation $\nu$ is uniquely determined.

Finally, we show that $\nu(x)\neq x$ for all $x\in V$ when $|V|>1$, which finishes the proof of the claim. Indeed, if $\nu(x)=x$ for some $x\in V$, then $V=I(x,\nu(x))=I(x,x)=\{x\}$, a contradiction to $|V|>1$.
\end{proof}

\begin{thm}\label{thm:self-inj}
$\Sigma_G$ is self-injective if and only if $G$ is interval-antipodal.
In such a case, the Nakayama permutation is involutive with no fixed point unless $|V|=1$.
\end{thm}

\begin{proof}
The algebra $\Sigma_G$ is self-injective if and only if every indecomposable projective module $P_x$ is injective. By Lemma~\ref{lem:proj-inj-max}, this is equivalent to the condition that for every $x\in V$ there exists a vertex $y\in V$ such that $I(x,y)=V$. By Lemma~\ref{lem:interval-antipodal-involution}, this is equivalent to $G$ being interval-antipodal.

Now assume that $G$ is interval-antipodal, and let $\nu:V\to V$ be the interval-antipodal permutation. Define a linear form $\lambda:\Sigma_G\to \Bbbk$ by
$$\lambda(b_{xy})=
\begin{cases}
1, & \text{if } y=\nu(x),\\
0, & \text{otherwise}.
\end{cases}$$
Since $I(x,\nu(x))=V$, for every $y\in V$ we have $b_{xy}b_{y,\nu(x)}=b_{x,\nu(x)}$. It follows that $\lambda$ is a Frobenius form on $\Sigma_G$. By abuse of notation, the Nakayama automorphism $\nu$ of $\Sigma_G$ with respect to $\lambda$ is given by
$\nu(b_{xy})=b_{\nu(x),\nu(y)}$ for all $x,y\in V$. Thus $\nu$ is precisely the automorphism induced by the interval-antipodal permutation $\nu$ on both indices. The last claim then follows from Lemma \ref{lem:interval-antipodal-involution}.
\end{proof}

\begin{prop}\label{prop:diam2 int-antipodal <=> cocktail party}
For a connected graph $G$ of diameter $2$, $G$ is interval-antipodal if and only if $G\cong K_{2n}\setminus M$ for a perfect matching $M$ on $K_{2n}$, i.e., a complete $n$-partite graph, a.k.a. the cocktail party graph.
In such a case, we have $d(v,\nu(v))=2$ and $\{v,\nu(v)\}\in M$ for all $v\in V(G)$.
In this case, $\Sigma_G$ is self-injective, and hence
$\operatorname{findim}\Sigma_G=\operatorname{idim}\Sigma_G=0$, $\operatorname{domdim}\Sigma_G=\infty$.
\end{prop}
\begin{proof}
We first show the claim that the antipode is at distance $2$.
Indeed, if $d(u,v)=1$, then $I(u,v)=\{u,v\}$, and so $I(u,v)=V(G)$ implies that $V(G)=\{u,v\}$, which contradicts the diameter $2$ assumption.  Hence, an interval-antipodal diameter $2$ graph must have $d(v,\nu(v))=2$ for all $v\in V(G)$.

$G$ being of diameter $2$ implies that $I(u,v)=\{u,v\}\cup (N(u)\cap N(v))$.  Hence, $d(u,v)=2$ with $I(u,v)=V(G)$ is equivalent to saying that $\{u,v\}$ is not an edge, and every other vertex in $G$ is adjacent to both $u$ and $v$.  In other words, the complement $K_{|V(G)|}\setminus E(G)$ is a disjoint union of single edges, which means that it is a perfect matching.
\end{proof}

\begin{rmk}
Since $K_{2n}\setminus M\cong K_{2,\ldots,2}$ is diagonal
\cite[Example~7.6]{HW17}, its distance algebra
\(\Sigma_G\) is Koszul. Moreover, $\Sigma_G$ is self-injective, and
hence its Koszul dual is an Artin--Schelter regular algebra of global dimension $2=\diam(G)$; see \cite[Proposition~5.10]{Smi96} and \cite{MV99}.
\end{rmk}

The following corollary follows directly from Theorem~\ref{thm:self-inj} and Proposition~\ref{prop:homological-dimensions-tiled-order}.

\begin{cor} 
\label{thm:tiled-gorenstein-antipodal}
For a finite connected graph $G$, the following conditions are
equivalent.
\begin{enumerate}
    \item $G$ is interval-antipodal.
    \item The tiled order $\Lambda_G$ is \(1\)-Iwanaga--Gorenstein.
    \item The distance algebra $\Sigma_G$ is self-injective.
\end{enumerate}
If these conditions hold and $\nu$ is the interval-antipodal
permutation, then $\nu$ is the Nakayama permutation of both
$\Lambda_G$ and $\Sigma_G$.
\end{cor}
\begin{rmk}
    In the classical use of terminology in the theory of orders, a \(1\)-Iwanaga--Gorenstein $R$-order in our sense is usually called a Gorenstein order; see, for example, \cite[Definition 6.1]{DKR67} and \cite[Section~2]{Rog73}.
    Note that the equivalence of (1) and (2) is already highlighted in \cite[Section~2.3]{IL26}, although only  stated in purely algebraic terms; in fact, the connection between $\Lambda_G$ and $\Sigma_G$ was first pointed out to us by Iyama.
\end{rmk}
\begin{rmk}
As explained in \cite{IL26}, for an interval-antipodal graph $G$, the Gorenstein parameter of the tiled order $\Lambda_G$ is always equal to $\diam(G)$; this is the order analogue of Lemma \ref{lem:Loewy legnth=ecc}.  The Gorenstein parameter being constant can also be deduced from a result of Martinez-Villa \cite[Theorem 3.3]{MV99}: for a graded self-injective algebra, the Loewy length is constant over all indecomposable projective modules.
\end{rmk}

We finish this subsection with the following remark.

\begin{prop}
The generalized Nakayama conjecture holds for distance algebras. More
precisely,
for every simple module $S$, there is some $i$ such that $\Ext_{\Sigma_G}^i(S,\Sigma_G)\neq 0$.
In particular, $\Sigma_G$ is self-injective if and only if it has infinite dominant dimension $\operatorname{domdim}\Sigma_G=\infty$.
\end{prop}

\begin{proof}
The path-length grading makes $\Sigma_G$ a finite-dimensional
non-negatively graded algebra with semisimple degree-zero part. 
Hence the generalized Nakayama conjecture for $\Sigma_G$ follows from \cite[Section~2]{Wil83}.

For completeness, we recall Wilson's argument in the following.
Let
$$\mathbb Z\{q\}:=\left\{
\sum_{r\geq r_0}a_rq^r
\;\middle|\;
r_0\in\mathbb Z,\ a_r\in\mathbb Z
\right\}=\mathbb Z[[q]][q^{-1}]$$
be the ring of Laurent series.
For a graded algebra $A$, the Grothendieck group $K_0(\mathrm{grmod} A)$ of the category of finitely generated graded $A$-modules is naturally a $\mathbb{Z}[q,q^{-1}]$-module, where $q$ acts as grading shift.

Consider the extension $\widehat{K}_A := \mathbb{Z}\{q\}\otimes_{\mathbb{Z}[q,q^{-1}]} K_0(\mathrm{grmod} A)$. Taking minimal graded projective resolutions of graded modules induces a $\mathbb{Z}\{q\}$-basis for $\widehat{K}_A$, given by the indecomposable projective $A$-modules generated in degree $0$.
Likewise, taking minimal graded injective coresolutions induces another basis given by the indecomposable injective $A$-modules cogenerated in degree $0$.
Thus, by a change-of-basis argument, every indecomposable injective module $I$ must occur in the minimal graded injective coresolution of $A$.  In other words, for every simple module $S$, we have $\Ext_{\mathrm{grmod}A}^i(q^jS,A)\neq 0$ for some $i,j\in \mathbb{Z}$.  Since $\Ext_A^k(S,A) = \bigoplus_{j\in\mathbb{Z}} \Ext_{\mathrm{grmod} A}^k(q^jS,A)$, the claim follows.   
\end{proof}

\subsection{The case of modular lattices}\label{subsec:modular-lattice}

We regard a finite poset $(P,\leq)$ as a graph by taking its undirected Hasse diagram. Thus $d(x,y)$ is defined for all $x,y\in P$, and we write $\Sigma_P$ for the distance algebra of this graph.

\begin{dfn}
Let $(P,\leq)$ be a finite poset. We say that $P$ is \defn{ranked} if
there exists a \defn{rank function}
$\rho:P\to\mathbb Z_{\geq0}$ such that, whenever $x$ covers $y$, one has
$\rho(x)=\rho(y)+1$.
\end{dfn}
We normalize the rank so that every minimal element has rank $0$. If
$x\leq y$ in a ranked poset, then $d(x,y)=\rho(y)-\rho(x)$.

Recall that a \defn{lattice} is a poset $L$ in which every pair
$x,y\in L$ has a least upper bound, called its \defn{join} and denoted
$x\vee y$, and a greatest lower bound, called its \defn{meet} and denoted
$x\wedge y$. We write $0$ and $1$ for the minimum and maximum of a finite lattice, respectively, except perhaps in the case of finite chains where we will often use specific labeling for each individual element.

\begin{dfn}
A lattice $L$ is
\begin{itemize}
    \item \defn{modular} if $x\leq y$ implies
    $x\vee(w\wedge y)=(x\vee w)\wedge y$ for any $w$;

    \item \defn{distributive} if meets are distributive over joins, i.e.,
    \[
    x\wedge(y\vee z)=(x\wedge y)\vee(x\wedge z).
    \]
    This is equivalent to joins being distributive over meets.
\end{itemize}
\end{dfn}
Note that modular lattices are always ranked. This is a classical result of Dedekind, which says that every maximal chain in an interval $[x,y]$ in such a lattice has the same length. Distributive lattices are modular. For example, a finite chain of length $m\ge 1$
\[
[m]=\{1<\cdots<m\}
\]
is clearly distributive.
Finite lattices (resp. modular lattices, resp. distributive lattices) are closed under direct product (same as the Cartesian product).  In particular, the Boolean lattice $[2]^n$ of rank $n\ge 0$ is distributive.
Note that if there is a poset isomorphism $P\cong P_1\times P_2$, then the Hasse diagram of $P$ is the Cartesian product of those of $P_1$ and $P_2$; in particular, we have $\Sigma_P \cong \Sigma_{P_1}\otimes_{\Bbbk} \Sigma_{P_2}$.

\begin{lem}\label{lem:modular-lattice-distance}
The following hold for any finite connected modular lattice.
\begin{enumerate}
    \item $d(x,y)=\rho(x)+\rho(y)-2\rho(x\wedge y)$;
    \item $d(x,y)=2\rho(x\vee y)-\rho(x)-\rho(y)$;
    \item $\rho(x)+\rho(y)=\rho(x\vee y)+\rho(x\wedge y)$;
    \item $d(x,y)=\rho(x\vee y)-\rho(x\wedge y)$.
\end{enumerate}
\end{lem}

\begin{proof}
(1) follows from $d(x,x\wedge y)=\rho(x)-\rho(x\wedge y)$. Part (2) is dual to (1), and (3) follows by combining (1) and (2). Substituting (3) into (1), or equivalently into (2), gives (4).
\end{proof}

Recall that in a finite lattice $L$, an element $x\in L$ is \defn{complemented} if there exists $y\in L$ such that \[ x\wedge y =0, \quad x\vee y=1.\]
In this case, $y$ is called the \defn{complement} of $x$.

\begin{lem}\label{lem:modular-lattice-central-antipode}
Let $L$ be a finite modular lattice and let $x\in L$. The following
conditions are equivalent.
\begin{enumerate}
\item $P_x \cong I_y$ as $\Sigma_L$-modules.
    \item There exists $y\in L$ such that $I(x,y)=L$.
    \item There exists $y\in L$ such that $x\wedge y=0$ and
    \[
    z=(x\wedge z)\vee(y\wedge z)
    \qquad\text{for every }z\in L.
    \]
    \item The element $x$ is \defn{central}\footnote{Properly, one should first define ``neutrality'' and then central means a neutral element with a complemented neutral element; see \cite[III.1]{Gra03} for details.}, that is, $x$ has a complement $y$ such that
    \[
    \Phi:[0,x]\times[0,y]\longrightarrow L,
    \qquad (u,v)\longmapsto u\vee v,
    \]
    is a lattice isomorphism.
\end{enumerate}
Consequently, $x$ has an antipode in $L$ if and only if $x$ is central. 
In that case, its antipode is its complement.
\end{lem}

\begin{proof}
$(1)\Leftrightarrow(2)$: This is Lemma \ref{lem:proj-inj-max}.

$(2)\Rightarrow(3)$: Assume that $I(x,y)=L$. Since $0\in I(x,y)$,
Lemma~\ref{lem:modular-lattice-distance}~(1) gives
\[
\rho(x)+\rho(y)-2\rho(x\wedge y)
=d(x,y)=d(x,0)+d(0,y)=\rho(x)+\rho(y).
\]
Thus $x\wedge y=0$. For any $z\in L$, the equality
$d(x,y)=d(x,z)+d(z,y)$ and the same distance formula give
$\rho(z)=\rho(x\wedge z)+\rho(y\wedge z)$. The meet of $x\wedge z$ and
$y\wedge z$ is $0$. By the rank identity in
Lemma~\ref{lem:modular-lattice-distance}~(3), their join has rank
$\rho(z)$. Since that join is at most $z$, it equals $z$.

$(3)\Rightarrow(2)$: For $z\in L$, put $u=x\wedge z$ and
$v=y\wedge z$. Then $u\wedge v=0$ and $z=u\vee v$.
Lemma~\ref{lem:modular-lattice-distance}~(3) therefore gives
$\rho(z)=\rho(u)+\rho(v)$. Using
Lemma~\ref{lem:modular-lattice-distance}~(1) and $x\wedge y=0$, we obtain
\[
d(x,z)+d(z,y)=\rho(x)+\rho(y)=d(x,y).
\]
Thus every $z\in L$ belongs to $I(x,y)$.

$(3)\Rightarrow(4)$: Define
\[
\Psi:L\longrightarrow[0,x]\times[0,y],
\qquad z\longmapsto(x\wedge z,y\wedge z).
\]
Condition (3) says that $\Phi\Psi=\operatorname{id}_L$. If $u\leq x$
and $v\leq y$, the modular law and $x\wedge y=0$ give
\[
x\wedge(u\vee v)=u,
\qquad
y\wedge(u\vee v)=v.
\]
Hence $\Psi\Phi$ is also the identity. Since both maps are
order-preserving, $\Phi$ is a lattice isomorphism.

$(4)\Rightarrow(3)$: Since $\Phi$ is a lattice isomorphism, projecting $z\in L$ onto a component $[0,w]$ for $w\in\{x,y\}$ is given by $z\wedge w$.  Hence, composing the projection with $\Phi$ yields (3).
\end{proof}

Lemma \ref{lem:modular-lattice-central-antipode} says that the self-injective part of $\Sigma_L$ for a modular lattice $L$ is controlled by its \defn{center} 
$$ Z(L):=\{ x \in L\mid \, x \text{ is a central element}\}$$
of $L$.  
Recall that a poset $P$ is \defn{directly indecomposable} if $P\cong P_1\times P_2$ implies one of $P_i$'s is a singleton.  Note that the direct product decomposition of a finite poset is unique \cite{Has51}.  

\begin{cor}\label{cor:directly indecomposable lattice}
The following are equivalent for a finite modular lattice $L$.
\begin{enumerate}
    \item $L$ is directly indecomposable.
    \item $Z(L)=\{0,1\}$.
    \item The only indecomposable projective-injective $\Sigma_L$-modules are $P_1\cong I_0$ and $P_0\cong I_1$.
\end{enumerate}
\end{cor}
\begin{proof}
$(1)\Leftrightarrow(2)$ is immediate by the definition of a central element and of direct indecomposability.  $(2)\Leftrightarrow(3)$ is a consequence of Lemma \ref{lem:modular-lattice-central-antipode} (1)$\Leftrightarrow$(4).
\end{proof}
\begin{rmk}\label{rem:Z(L)=Boolean}
Thanks to the unique decomposition property, by writing a finite modular lattice $L$ as a direct product of $n$ directly indecomposable lattices $L_1,\ldots, L_n$, we have $$Z(L) \cong Z(L_1\times \cdots  \times  L_n) \cong Z(L_1)\times \cdots \times Z(L_n) \cong [2]^n.$$  
In particular, the \defn{atoms} (elements covering $0$) of $Z(L)\cong [2]^n$ are in one-to-one correspondence with the directly indecomposable factors of $L$, and the complement of each atom is a \defn{coatom} (element covered by $1$).    
\end{rmk}

\begin{cor}\label{cor:modular-lattice-interval-antipodal}
The following are equivalent for a finite modular lattice $L$.
\begin{enumerate}
    \item $\Sigma_L$ is self-injective.
    \item $L$ is interval-antipodal.
    \item $L\cong Z(L)$.
    \item $L$ is a Boolean lattice.
\end{enumerate}
\end{cor}

\begin{proof}
(1) $\Leftrightarrow$ (2) is Theorem \ref{thm:self-inj}. 
(2) $\Leftrightarrow$ (3) follows from the equivalence between (1) and (4) of Lemma \ref{lem:modular-lattice-central-antipode}.  (3) $\Leftrightarrow$ (4) follows from Remark \ref{rem:Z(L)=Boolean}.
\end{proof}

Recall that a graph $G$ is called \defn{modular} if $I(u,v)\cap I(v,w)\cap I(w,u)\neq\emptyset$ for every three vertices $u,v,w$. The Hasse diagram of a modular lattice is a standard example of a modular graph. However, an interval-antipodal modular graph need not be a hypercube (Hasse diagram of a Boolean lattice).

\begin{eg}\label{eg:non-cubical-interval-antipodal}
Let $K_{5\times 2}:=K_{5,5}\setminus M$, where the bipartition is given by $\{x_1,\ldots,x_5\}$ and $\{y_1,\ldots,y_5\}$, and $M=\{x_i y_i\mid 1\leq i\leq5\}$. Thus $x_i$ is adjacent to $y_j$ precisely when $i\neq j$.

\begin{center}
\begin{tikzpicture}[
  x=1cm,
  y=1cm,
  vertex/.style={
    circle,
    draw,
    fill=white,
    inner sep=0pt,
    minimum size=7mm
  },
  edge/.style={draw=black,line width=.4pt}
]
  \coordinate (x1) at (-4, 1.4);
  \coordinate (x2) at (-2, 1.4);
  \coordinate (x3) at ( 0, 1.4);
  \coordinate (x4) at ( 2, 1.4);
  \coordinate (x5) at ( 4, 1.4);

  \coordinate (y1) at (-4,-1.4);
  \coordinate (y2) at (-2,-1.4);
  \coordinate (y3) at ( 0,-1.4);
  \coordinate (y4) at ( 2,-1.4);
  \coordinate (y5) at ( 4,-1.4);

  \draw[edge] (x1)--(y2);
  \draw[edge] (x1)--(y3);
  \draw[edge] (x1)--(y4);
  \draw[edge] (x1)--(y5);

  \draw[edge] (x2)--(y1);
  \draw[edge] (x2)--(y3);
  \draw[edge] (x2)--(y4);
  \draw[edge] (x2)--(y5);

  \draw[edge] (x3)--(y1);
  \draw[edge] (x3)--(y2);
  \draw[edge] (x3)--(y4);
  \draw[edge] (x3)--(y5);

  \draw[edge] (x4)--(y1);
  \draw[edge] (x4)--(y2);
  \draw[edge] (x4)--(y3);
  \draw[edge] (x4)--(y5);

  \draw[edge] (x5)--(y1);
  \draw[edge] (x5)--(y2);
  \draw[edge] (x5)--(y3);
  \draw[edge] (x5)--(y4);

  \node[vertex] at (x1) {$x_1$};
  \node[vertex] at (x2) {$x_2$};
  \node[vertex] at (x3) {$x_3$};
  \node[vertex] at (x4) {$x_4$};
  \node[vertex] at (x5) {$x_5$};

  \node[vertex] at (y1) {$y_1$};
  \node[vertex] at (y2) {$y_2$};
  \node[vertex] at (y3) {$y_3$};
  \node[vertex] at (y4) {$y_4$};
  \node[vertex] at (y5) {$y_5$};
\end{tikzpicture}
\end{center}
Then $K_{5\times 2}$ is clearly not a hypercube since $10\neq 2^n$ for any $n$.

Note that for every $i$ and every vertex $z$, we have
$d(x_i,z)+d(z,y_i)=3=d(x_i,y_i)$. Thus $I(x_i,y_i)=V(G)$, and $y_i$ is the unique farthest vertex from $x_i$. Hence $G$ is interval-antipodal, with antipodes $x_i\leftrightarrow y_i$.
It remains to check that $K_{5\times 2}$ is modular. If two of the three vertices
coincide, their three pairwise intervals clearly have a common
vertex, so suppose that the vertices are distinct. If they are
$x_i,x_j,x_k$, choose $\ell\notin\{i,j,k\}$. The paths
$x_i-y_\ell-x_j$, $x_j-y_\ell-x_k$, and
$x_k-y_\ell-x_i$ are shortest paths. Hence
\[
y_\ell\in
I(x_i,x_j)\cap I(x_j,x_k)\cap I(x_k,x_i).
\]
The case of three vertices of the form $y_i,y_j,y_k$ is the same. It remains to consider $x_i,x_j,y_k$ with $i\neq j$; the case of one vertex of the form $x_i$ and two of the form $y_j$ is analogous. If $k\notin\{i,j\}$, then $x_i-y_k-x_j$ is a shortest path, so
\[
y_k\in
I(x_i,x_j)\cap I(x_j,y_k)\cap I(y_k,x_i).
\]
If $k=i$, choose $\ell\notin\{i,j\}$. Then
$x_i-y_\ell-x_j-y_i$ is a shortest path, and therefore
\[
x_j\in
I(x_i,x_j)\cap I(x_j,y_i)\cap I(y_i,x_i).
\]
The case $k=j$ follows in the same way, with $x_i$ in the
intersection. Thus $K_{5\times 2}$ is modular.
\end{eg}

\section{Homological dimensions}\label{sec:tree}

\subsection{General elementary facts}

Following \cite[Definition~2.1(4) and Remark~2.2]{AW24}, a graph $G$ is called \defn{geodetic} if any two vertices are joined by a unique shortest path. Note that the distance algebra $\Sigma_G$ is a monomial algebra if and only if $G$ is geodetic.

\begin{lem}\label{lem:mono-findim}
Assume that the distance algebra $\Sigma_G$ is monomial. Then $\operatorname{Findim}\Sigma_G\leq1$.
In particular, $\operatorname{findim}\Sigma_G\leq1$.
\end{lem}

\begin{proof}
Let $M$ be an arbitrary right $\Sigma_G$-module with finite
projective dimension. Since $\Sigma_G$ is a monomial algebra, we use the standard fact that the second syzygy $\Omega^2(M)$ is a direct sum of right ideals generated by non-trivial paths; see \cite[Theorem I]{HZ91}. Hence
$$\Omega^2(M)\cong \bigoplus_{\lambda} p_\lambda \Sigma_G$$
for some non-zero non-trivial paths $p_\lambda$.

In the distance algebra, every non-zero path is represented by some basis element $b_{xy}$. Since the paths $p_\lambda$ are non-trivial, we have $x\neq y$. By Proposition \ref{prop:pdim Pxy}, every $P_x^y=b_{xy}\Sigma_G$ with $x\neq y$ has infinite projective dimension. Since $\Omega^2(M)$ has finite projective dimension, no such summand can occur. Therefore $\Omega^2(M)=0$. Thus $\operatorname{pdim}_{\Sigma_G}M\leq 1$.
Hence every module of finite projective dimension has projective dimension at most $1$, and consequently $\operatorname{findim} \Sigma_G\leq 1$.
\end{proof}
\begin{rmk}
The description of $\Omega^2(M)$ above is closely related to the study in \cite[Section~2]{GHZ98}; see \cite[Example~2.3(a)]{GHZ98}. Moreover, if the tiled order $\Lambda_G$ has finite global dimension, then every finitely generated $\Sigma_G\cong\Lambda_G/x\Lambda_G$-module $M$ has finite syzygy type; that is, up to isomorphism, only finitely many indecomposable modules occur as direct summands of the syzygies of $M$; see \cite[Corollary~4.4]{GHZ98}.
\end{rmk}

We recall the following elementary facts.
\begin{prop}\label{prop:homological-dimensions-tensor}
Let $A$ and $B$ be finite-dimensional $\Bbbk$-algebras. Then the following hold.
\begin{enumerate}
    \item {\rm 
\cite[Theorem~16]{ERZ57}} $\operatorname{findim}(A\otimes_{\Bbbk}B)=\operatorname{findim}A+\operatorname{findim}B$.
\item {\rm\cite[Theorem~XI.3.1]{CE56}} $\operatorname{idim}(A\otimes_{\Bbbk}B)=\operatorname{idim}A+\operatorname{idim}B$.
\item {\rm \cite[Lemma~6]{Mul68}} $\operatorname{domdim}(A\otimes_{\Bbbk}B)=\min\{\operatorname{domdim}A,\operatorname{domdim}B\}$.
\end{enumerate}
Here the equalities are understood in $\mathbb Z_{\geq 0}\cup\{\infty\}$.
\end{prop}

\begin{rmk}
    By Propositions~\ref{prop:tensor} and \ref{prop:homological-dimensions-tensor}, Cartesian products of graphs yield distance algebras with arbitrarily large finite finitistic and injective dimensions. In contrast, taking Cartesian products cannot increase the dominant dimension.
\end{rmk}

We now turn to the computation of these homological dimensions for several specific families of graphs.

\subsection{The case of tree graphs}\label{subsec:tree}

We recall the definition of a preprojective algebra. Let $Q$ be a
finite quiver, and let $\overline Q$ be its double quiver, obtained
from $Q$ by adjoining, for every arrow $a:i\to j$, an arrow
$a^*:j\to i$. The \defn{preprojective algebra} of $Q$ is
$$\Pi(Q):=\Bbbk\overline Q\bigg/\left\langle\sum_{a\in Q_1}(aa^*-a^*a)
\right\rangle.$$
Up to isomorphism, $\Pi(Q)$ depends only on the underlying graph of
$Q$ and not on the chosen orientation. Preprojective algebras were
introduced by Gelfand and Ponomarev in \cite{GP79}.

\begin{prop}
If $T$ is a tree, then $\Sigma_T$ is a quotient of $\Pi(T)$.
\end{prop}

\begin{proof}
Let $Q_T$ be the double quiver of $T$. Note that the preprojective algebra $\Pi(T)$ of $T$ is $\Bbbk Q_T/\langle \rho\rangle$ where
$$\rho=\sum_{\{x,y\}\in E(T)}a_{xy}a_{yx}-a_{yx}a_{xy}.$$
Hence, $\Sigma_T \cong \Pi(T)/I$ where $I=\langle a_{xy}a_{yx}\mid a_{xy}\in Q_T\rangle$.
\end{proof}

While preprojective algebras are quite well-studied, this particular quotient is not, but there is some use of it to obtain results for preprojective algebras, for instance, in \cite{Miz14}.

\begin{prop}\label{prop:idim Sigma tree}
Let $T$ be a tree. Then
$$
\operatorname{Findim} \Sigma_T=\operatorname{findim} \Sigma_T=\operatorname{idim}\Sigma_T=
\begin{cases}
0, & |V(T)|\leq 2,\\
1, & |V(T)|\geq 3.
\end{cases}
$$
\end{prop}

\begin{proof}
If $|V(T)|\leq 2$, then $\Sigma_T$ is self-injective. Now assume $|V(T)|\geq 3$. Since $|V(T)|\geq 3$, the tree $T$ is not interval-antipodal. Equivalently, by Theorem \ref{thm:self-inj}, $\Sigma_T$ is not self-injective. Hence
$\operatorname{idim}\Sigma_T\neq 0$. We just need to show that $\operatorname{idim}\Sigma_T\leq 1$. 

Now we fix a vertex $x\in V(T)$. We regard $T$ as a rooted tree with root $x$. For a vertex $v\in V(T)$, set
$$\operatorname{ch}_x(v):=\{w\in V(T)\mid w \text{ is adjacent to } v\text{ and } d(x,w)=d(x,v)+1\}.$$
Thus $\operatorname{ch}_x(v)$ is the set of children of $v$ with respect to the root $x$. Let $c_x(v):=|\operatorname{ch}_x(v)|$.

Let $I_x$ be the injective envelope of the simple module $S_x$, that is, $I_x=D(\Sigma_Te_x)$. We claim that every $I_x$ has projective dimension at most one.  For each vertex $v\in V(T)$, let $\xi_v\in I_x$ be the dual basis element corresponding to $b_{vx}\in \Sigma_Te_x$. Thus $I_x=\operatorname{span}_{\Bbbk}\{\xi_v\mid v\in V(T)\}$. For each $y\in \Max(x)$, define a homomorphism $\pi_y:e_y\Sigma_T\to I_x$ induced by $\pi_y(e_y)=\xi_y$. Now set
$$\pi:=\bigoplus_{y\in \Max(x)}\pi_y:\bigoplus_{y\in \Max(x)}e_y\Sigma_T\to I_x.$$
This map is surjective. Indeed, every vertex $v\in V(T)$ lies on some path from $x$ to some vertex in $\Max(x)$, and hence $\xi_v$ lies in the image of some $\pi_y$.

It remains to describe $\Ker\pi$. Let $v\in V(T)$ with
$$\operatorname{ch}_x(v)=\{w_1,\dots,w_s\},\qquad s=c_x(v)\geq 2.$$
For each child $w_a$, choose a vertex $y_a\in \Max(x)$ lying in the branch starting at $w_a$. Such a vertex exists because each branch of a finite rooted tree contains a maximal vertex for the distance function $d(x,-)$. For $a=1,\dots,s-1$, define a homomorphism
$$\iota_{v,a}:e_v\Sigma_T\to \bigoplus_{y\in \Max(x)}e_y\Sigma_T$$
by sending the generator $e_v$ to $b_{y_av}-b_{y_sv}$, where $b_{y_av}$ is placed in the summand $e_{y_a}\Sigma_T$ and $b_{y_sv}$ is placed in the summand $e_{y_s}\Sigma_T$. Since both $y_a$ and $y_s$ lie in different branches starting at children of $v$, the paths from $y_a$ to $x$ and from $y_s$ to $x$ both pass through $v$. Therefore
$$\pi(b_{y_av})=\xi_v=\pi(b_{y_sv}),$$
and hence $\operatorname{Im}\iota_{v,a}\subseteq \Ker\pi.$

Taking the direct sum of these maps over all vertices $v$ with $c_x(v)\geq 2$, we obtain a homomorphism
$$\iota:\bigoplus_{v\in V(T)} (e_v\Sigma_T)^{\oplus \max\{c_x(v)-1,0\}}\to
\bigoplus_{y\in \Max(x)}e_y\Sigma_T.$$
By construction, $\pi\iota=0$. A basis check shows that $\operatorname{Im}\iota=\Ker\pi$. Indeed, the only relations among the elements mapping to a fixed basis vector $\xi_v$ come from the different branches above $v$; if $v$ has $c_x(v)$ children, these give exactly $c_x(v)-1$ independent relations, represented by the differences $b_{y_av}-b_{y_sv}$ with $a=1,\dots,c_x(v)-1$.
Thus we obtain the exact sequence that is the projective presentation of $I_x$
$$0\to\bigoplus_{v\in V(T)} (e_v\Sigma_T)^{\oplus \max\{c_x(v)-1,0\}}
\xrightarrow{\ \iota\ }\bigoplus_{y\in \Max(x)} e_y\Sigma_T\xrightarrow{\ \pi\ }I_x\to 0.$$
Hence $\operatorname{pdim}I_x\leq 1$ for all $x\in V(T)$. 

Now use the duality $\omega^{-1}D$ in Proposition \ref{prop:dual-functor}. Since $\omega^{-1}D(e_x\Sigma_T)\cong I_x$, we have
$$ \operatorname{idim}\Sigma_T=\sup_{x\in V(T)}\operatorname{idim}(e_x\Sigma_T)=\sup_{x\in V(T)}\operatorname{pdim}I_x\leq 1.$$
Therefore $\operatorname{idim}\Sigma_T=1$.

Therefore, $\Sigma_T$ is Iwanaga--Gorenstein, and hence
$\operatorname{Findim} \Sigma_T=\operatorname{findim} \Sigma_T=\operatorname{idim}\Sigma_T$.
\end{proof}

\begin{prop}\label{prop:domdim Sigma tree}
Let $T$ be a tree. Then
$$
\operatorname{domdim} \Sigma_T=
\begin{cases}
\infty, & |V(T)|\leq 2,\\
1, & T=A_n\text{ with }n\geq 3,\\
0, & \text{otherwise}.
\end{cases}$$
\end{prop}

\begin{proof}
If $|V(T)|\leq 2$, then $\Sigma_T$ is self-injective. Now assume $|V(T)|\geq 3$. Let $L$ be the set of leaves of $T$. For a tree, the local maxima of the distance function are exactly the leaves away from the starting vertex:
$$\Max(x)=
\begin{cases}
L, & x\notin L,\\
L\setminus\{x\}, & x\in L.
\end{cases}$$
By the computation of socles,
$\operatorname{soc} (e_x\Sigma_T)\cong\bigoplus_{y\in \Max(x)}S_y$.
Therefore the injective envelope of $e_x\Sigma_T$ is $I(e_x\Sigma_T)\cong\bigoplus_{y\in \Max(x)}I_y$.

Suppose first that $T=A_n$ with $n\geq 3$. Then $T$ has exactly two leaves, say $l$ and $r$. We have
$$I_l\cong e_r\Sigma_T,\qquad I_r\cong e_l\Sigma_T.$$
Thus $I_l$ and $I_r$ are projective. Since every $\Max(x)$ is contained in $\{l,r\}$, the injective envelope $I(\Sigma_T)$ is projective. Hence $\operatorname{domdim}\Sigma_T\geq 1$. However, for an internal vertex $x$, the injective module $I_x$ is not projective. Moreover, the cokernel of the injective envelope map $e_x\Sigma_T\rightarrow I(e_x\Sigma_T)\cong I_l\oplus I_r$ is isomorphic to $I_x$. Thus the next term in the minimal injective resolution of $\Sigma_T$ contains a non-projective injective direct summand. Therefore $\operatorname{domdim}\Sigma_T=1$.

Finally, suppose that $T$ is not a path graph. Then $T$ has at least three leaves. For any leaf $l\in L$, the injective module $I_l$ is not projective. 
Since every $\Max(x)$ contains leaves, the injective envelope of the regular module contains some non-projective injective module $I_l$ as a direct summand. Hence $I(\Sigma_T)$ is not projective. Therefore $\operatorname{domdim}\Sigma_T=0$.
\end{proof}

By Theorem~\ref{thm:rep-finite}, the representation-finite distance algebras are precisely the distance algebras associated with path graphs. Moreover, for $n\geq3$, the preceding results show that $\Sigma_{A_n}$ is $0$-minimal Auslander--Gorenstein. 
Now we can determine the dominant dimension of \emph{every} indecomposable module over $\Sigma_{A_n}$. 
Recall that we have the modules
\[
P_{i+1}^{\,i}:=b_{i+1,i}\Sigma_{A_n}
\quad(1\leq i<n),
\qquad
P_{i-1}^{\,i}:=b_{i-1,i}\Sigma_{A_n}
\quad(2\leq i\leq n).
\]
In the following, we write $M\in \mathcal{X}$ if $M$ is isomorphic to any one of the indecomposable modules appearing in $\mathcal{X}$.

\begin{prop}\label{prop:domdim-modules-An}
Let $n\geq3$. For every indecomposable $\Sigma_{A_n}$-module $M$, one has
\[
\operatorname{domdim}M=
\begin{cases}
\infty, & M\in\{ P_1, P_n\},\\
2,      & M\in\{ P_1^{\,2}, P_n^{\,n-1} \}, \\
1,      & M\in \{ P_{i+1}^{\,i}, P_j, P_{k-1}^k \mid 1\le i<n-1, 1<j<n, 2<k\le n \}\\
0,      & \text{otherwise}.
\end{cases}
\]
In particular, there are precisely $3n-2$ isomorphism classes of
indecomposable $\Sigma_{A_n}$-modules of positive dominant dimension.
\end{prop}

\begin{proof}
By Example~\ref{eg:line and cycle}(1), the algebra $\Sigma_{A_n}$ is
gentle, and by Proposition~\ref{prop:biserial rep type} it has no band modules. Hence every indecomposable $\Sigma_{A_n}$-module is a string module. Consequently, every indecomposable module
is a thin interval module supported on an interval $[l,r]\subseteq[1,n]$, whose quiver induced by its string is an orientation of that interval.

Let $\operatorname{Sink}(M)$ denote the set of sinks in the quiver of $M$. Then
$\operatorname{soc}M
 \simeq\bigoplus_{j\in\operatorname{Sink}(M)}S_j$
and therefore $I(M)\simeq
\bigoplus_{j\in\operatorname{Sink}(M)}I_j.$
A direct calculation gives $I_1\simeq P_n$ and $I_n\simeq P_1$.
Thus $I_1$ and $I_n$ are projective-injective. On the other hand,
$I_j$ is not projective for $1<j<n$. It follows that $I(M)$ is projective if and only if $\operatorname{Sink}(M)\subseteq\{1,n\}$.

We first classify the oriented intervals satisfying this condition.
Every finite oriented interval has a sink, so its support must contain at least one endpoint of $A_n$. Suppose that the support is $[1,r]$ with $r<n$. The absence of an internal sink forces all arrows to point towards $1$, and the corresponding module is $P_{r+1}^{\,r}$. Dually, if the support is $[l,n]$ with $l>1$, all arrows must point towards $n$, and the corresponding module is $P_{l-1}^{\,l}$.

It remains to consider modules supported on $[1,n]$. Such an
orientation has no internal sink if and only if its direction sequence contains no pattern $(\longrightarrow\ \longleftarrow)$.
Hence, for a unique $i$, its coefficient quiver has the form
\[
1\longleftarrow\cdots\longleftarrow i
 \longrightarrow\cdots\longrightarrow n.
\]
This is precisely the indecomposable projective module $P_i$.
Therefore, every remaining indecomposable module has dominant
dimension zero.

It remains to compute the dominant dimensions of the modules just
listed. 

For $1<i<n$, the minimal injective coresolution begins
\[
0\longrightarrow P_i
 \longrightarrow I_1\oplus I_n
 \longrightarrow I_i
 \longrightarrow0.
\]
The first injective term is projective, whereas $I_i$ is not.
Consequently, $\operatorname{domdim}P_i=1$.

For $1\leq i<n$, there is an injective-envelope sequence
\[
0\longrightarrow P_{i+1}^{\,i}
 \longrightarrow I_1
 \longrightarrow C_i
 \longrightarrow0,
\]
where $C_i$ is the thin interval module supported on $[i+1,n]$ and
$\operatorname{soc}C_i\simeq S_{i+1}$. If $i<n-1$, then $I(C_i)=I_{i+1}$ is not projective, and hence $\operatorname{domdim}P_{i+1}^{\,i}=1$. For $i=n-1$, one has $C_{n-1}\simeq S_n$, and the minimal injective coresolution begins
\[
0\longrightarrow P_n^{\,n-1}
 \longrightarrow I_1
 \longrightarrow I_n
 \longrightarrow I_{n-1}
 \longrightarrow\cdots.
\]
The first two injective terms are projective, while $I_{n-1}$ is not. Thus $\operatorname{domdim}P_n^{\,n-1}=2$.
The dual argument gives
\[
\operatorname{domdim}P_{i-1}^{\,i}=1
\quad(2<i\leq n),
\qquad
\operatorname{domdim}P_1^{\,2}=2.
\]
This proves the asserted classification.
\end{proof}

\begin{rmk}\label{rmk:torsionfree-class-An}
By Proposition~\ref{prop:domdim-modules-An} and the (Auslander--Gorenstein version of) Auslander-type correspondence \cite[Theorem 4.4.1]{Iya07}, the full subcategories
\begin{align*}
\mathcal T_n
&:=\{M\in\operatorname{mod}\Sigma_{A_n}
  \mid Me_1=Me_n=0\}\simeq \operatorname{mod} \Sigma_{A_{n-2}}\\
\text{and}\quad  
\mathcal F_n &:=\{M\in\operatorname{mod}\Sigma_{A_n}
  \mid \operatorname{domdim}M\geq1\}
\end{align*}
form a hereditary torsion pair $(\mathcal T_n,\mathcal F_n)$ with $\mathcal F_n$ faithful. 
\end{rmk}

\subsection{The case of self-centered graphs}\label{sec:self-centered}

\begin{dfn}\label{def:self-centered}\cite{Buc89,BH90}
A graph $G$ is called \defn{self-centered} if $\operatorname{ecc}(x)=\operatorname{diam}(G)$, for all $x\in V(G)$, i.e., all vertices have the same eccentricity.
\end{dfn}

\begin{eg}\label{eg:self-centered eg}
Cycles $C_n$, complete graphs $K_n$, and complete bipartite graphs $K_{m,n}$ ($m,n\geq 2$) are all self-centered. 
\end{eg}

\begin{prop}
Every interval-antipodal graph is self-centered.
\end{prop}
\begin{proof}
Let $G$ be interval-antipodal and fix $x\in V(G)$. Since $I(x,\nu(x))=V(G)$,
for any $y,z\in V(G)$ we have $d(y,\nu(x))=d(x,\nu(x))-d(x,y)$ and $d(z,\nu(x))=d(x,\nu(x))-d(x,z)$. Therefore,
\begin{align*}
d(y,z)
&\leq \min\bigl\{
d(y,x)+d(x,z),\,
d(y,\nu(x))+d(\nu(x),z)
\bigr\}\\
&=\min\bigl\{
d(x,y)+d(x,z),\,
2d(x,\nu(x))-d(x,y)-d(x,z)
\bigr\}\\
&\leq d(x,\nu(x)).
\end{align*}
Thus $\operatorname{diam}(G)\leq d(x,\nu(x))$. The reverse inequality is immediate, and hence $d(x,\nu(x))=\operatorname{diam}(G)$. It follows that $\operatorname{ecc}(x)=\operatorname{diam}(G)$. Since $x$ was arbitrary, $G$ is self-centered.
\end{proof}

We record a simple consequence of the join construction. Recall from Section~\ref{sec:biserial} that $\overline N_G(x)$ denotes the set of non-neighbors of $x$ in $G$.

\begin{lem}\label{lem:join-eccentricity}
Let $K=G\vee H$, where $V(G)$ and $V(H)$ are non-empty. Then, for
$x\in V(G)$,
$$\operatorname{ecc}_K(x)=
\begin{cases}
1, & \text{if } \overline N_G(x)=\emptyset,\\
2, & \text{if } \overline N_G(x)\neq\emptyset.
\end{cases}
$$
\end{lem}

\begin{proof}
By Lemma~\ref{lem:join-distance}, every vertex in $V(G)$ is at distance
one from every vertex in $V(H)$. Moreover, two distinct vertices in $V(G)$ are at distance one in $K$ if they are adjacent in $G$, and at distance two otherwise. Hence the largest distance from $x\in V(G)$ is one precisely when $x$ is adjacent to every other vertex of $G$, and is two otherwise. The proof for vertices of $H$ is identical.
\end{proof}

\begin{cor}\label{cor:join-self-centered}
Let $K=G\vee H$, where $V(G)$ and $V(H)$ are non-empty. Then $K$ is
self-centered if and only if one of the following holds:
\begin{enumerate}
    \item $G$ and $H$ are both complete graphs, that is, $\overline N_G(x)=\emptyset$ and $\overline N_H(y)=\emptyset$ for all $x\in V(G)$, $y\in V(H)$. In this case, the join is just a complete graph.
    \item Every vertex of $G$ has a non-neighbor in $G$, and every
    vertex of $H$ has a non-neighbor in $H$, that is, $\overline N_G(x)\neq\emptyset$ and $\overline N_H(y)\neq\emptyset$ for all $x\in V(G)$, $y\in V(H)$.
\end{enumerate}
\end{cor}

\begin{proof}
By Lemma~\ref{lem:join-eccentricity}, every vertex of $K$ has eccentricity either one or two. Hence $K$ is self-centered exactly when all vertices have eccentricity one, or all vertices have eccentricity two. The first case means that $G$ and $H$ are both complete. The second case means that every vertex has a non-neighbor inside its own part.
\end{proof}

\begin{eg}
Take $G=H=C_4$. Every vertex then has a non-neighbor in its own copy of $C_4$, and hence every vertex of $C_4\vee C_4$ has eccentricity $2$. Thus $C_4\vee C_4$ is self-centered. Moreover, it is interval-antipodal.
\end{eg}
    
\begin{lem}\label{lem:join-interval-antipodal}
Let $K=G\vee H$, where $V(G)$ and $V(H)$ are non-empty. Assume
$|V(K)|>2$. Then $K$ is interval-antipodal if and only if $$
|\overline N_G(x)|=1\quad\text{for all }x\in V(G), \qquad|\overline N_H(y)|=1\quad\text{for all }y\in V(H).$$
\end{lem}

\begin{proof}
Let $x\in V(G)$. Since $|V(K)|>2$, if $z$ is adjacent to $x$ in $K$,
then $I_K(x,z)\neq V(K)$, because $d_K(x,z)=1$ and the interval between
two adjacent vertices contains only these two vertices. Thus an antipode of $x$, if it exists, must be a vertex $x'\in V(G)$ that is not adjacent to $x$ in $G$. Now suppose $x'\in \overline N_G(x)$. Then $d_K(x,x')=2$. For every $y\in V(H)$, we have
$$d_K(x,y)+d_K(y,x')=1+1=2,$$
so $V(H)\subseteq I_K(x,x')$. For a vertex $u\in V(G)\setminus\{x,x'\}$,
we have $u\in I_K(x,x')$ if and only if $d_K(x,u)=1$ and $d_K(u,x')=1$, that is, if and only if $u$ is adjacent in $G$ to both $x$ and $x'$. Therefore $I_K(x,x')=V(K)$ if and only if $x'$ is the unique non-neighbor of $x$ in $G$, and $x$ is the unique non-neighbor of $x'$ in $G$. It follows that every vertex of $V(G)$ has an antipode if and only if every vertex of $G$ has exactly one non-neighbor in $G$.
\end{proof}

\begin{eg}\label{eg:self-centered-revisited}
By Lemma \ref{lem:join-interval-antipodal}, the complete bipartite graph $K_{m,n}$ is interval-antipodal if and only if $m=n=2$. In particular, $K_{2,n}$ is self-centered but not interval-antipodal for every $n\geq3$.
\end{eg}

Thus, joins of graphs provide us with many examples of self-centered graphs. Next, we compute the homological dimensions of the distance algebras associated with self-centered graphs.

\begin{prop}
Let $G$ be a self-centered graph. Then
\[
\operatorname{Findim} \Sigma_G=\operatorname{findim} \Sigma_G=0.
\]
\end{prop}

\begin{proof}
Set $J=\operatorname{rad}\Sigma_G$, and $D=\operatorname{diam}(G)$.
Thus $J^r=\operatorname{span}_{\Bbbk}\{b_{xy}\mid d(x,y)\ge r\}$. Since $G$ is self-centered, every vertex $x$ satisfies $\operatorname{ecc}(x)=D$. Consequently, $J^{D+1}=0$ and, for every non-zero projective right $\Sigma_G$-module $P$, $PJ^D\neq 0$. Suppose that an arbitrary right $\Sigma_G$-module $M$ has finite projective dimension $r\geq 1$. Consider a minimal projective resolution
\[
0\longrightarrow P_r
 \xrightarrow{d_r}P_{r-1}
 \longrightarrow\cdots\longrightarrow P_0
 \longrightarrow M\longrightarrow 0.
\]
The map $d_r$ is injective, and minimality implies that $d_r(P_r)\subseteq P_{r-1}J$. Therefore
$$d_r(P_rJ^D)=d_r(P_r)J^D\subseteq P_{r-1}J^{D+1}=0.$$
However, $P_rJ^D\neq0$, contradicting the injectivity of $d_r$. Thus every finitely generated module of finite projective dimension is projective, and hence $\operatorname{findim}\Sigma_G=0$.
\end{proof}

\begin{prop}\label{prop:idim self-centered}
Let $G$ be a self-centered graph. Then
\[
\operatorname{idim} \Sigma_G=
\begin{cases}
0,      & \text{if $G$ is interval-antipodal},\\
\infty, & \text{if $G$ is not interval-antipodal}.
\end{cases}
\]
\end{prop}

\begin{proof}
By Theorem \ref{thm:self-inj}, if $G$ is interval-antipodal, then $\Sigma_G$ is self-injective and $\operatorname{idim}\Sigma_G=0$. 
Assume now that $G$ is not interval-antipodal. Suppose, for a  contradiction, that $\operatorname{idim}\Sigma_G<\infty$. Then by the duality in Proposition \ref{prop:dual-functor}, every indecomposable injective module $I$ has $\operatorname{pdim}I<\infty$. By the previous proposition, $\operatorname{findim}\Sigma_G=0$. Thus every such injective module must be projective. It follows that $\Sigma_G$ is
self-injective, contradicting the assumption that $G$ is not
interval-antipodal. Hence $\operatorname{idim}\Sigma_G=\infty$.
\end{proof}

\begin{eg}\label{eg:idim K2,3}
Let $G=K_{2,3}$ with the bipartition given by $\{0,1\}\sqcup \{a,b,c\}$.
Consider the indecomposable injective module $I_a$. The syzygy $\Omega(I_a)$ is a $5$-dimensional indecomposable module of Loewy length $2$ with top $S_0\oplus S_1$ and socle $S_{a}\oplus S_{b} \oplus S_{c}$.
For the second syzygy, we have $\Omega^2(I_a) \cong \omega_*^{-1}D(\Omega(I_a))$.  A direct (careful) calculation shows that $\Omega^3(I_a) \cong \Omega(I_a)^{\oplus 2}$, and so $\pdim I_a = \infty$.  In fact, we have
\[
\Omega^{2k+1}(I_a) \cong \Omega(I_a)^{\oplus 2^k} \text{ and }\Omega^{2k}(I_a) \cong \Omega^{2}(I_a)^{\oplus 2^{k-1}}
\]
for all $k\ge 1$.
\end{eg}

\begin{prop}
Let $G$ be a self-centered graph. Then
\[
\operatorname{domdim} \Sigma_G=
\begin{cases}
\infty, & \text{if $G$ is interval-antipodal},\\
0,      & \text{if $G$ is not interval-antipodal}.
\end{cases}
\]
\end{prop}

\begin{proof}
By Theorem \ref{thm:self-inj}, if $G$ is interval-antipodal, then $\Sigma_G$ is self-injective and $\operatorname{domdim}\Sigma_G=\infty$.

Suppose that $G$ is not interval-antipodal. We first observe that every simple right $\Sigma_G$-module occurs in the socle of the regular module. Indeed, let $y\in V(G)$. Since $G$ is self-centered, there exists
$x\in V(G)$ such that $d(x,y)=\operatorname{diam}(G)$. Thus $y\in\Max(x)$. Since
\[\operatorname{soc}(e_x\Sigma_G)\cong \bigoplus_{z\in\Max(x)}S_z,\]
the simple module $S_y$ occurs in $\operatorname{soc}\Sigma_G$. It follows that the injective envelope $I(\Sigma_G)$ contains every indecomposable injective right $\Sigma_G$-module as a direct summand.

Suppose that $\operatorname{domdim}\Sigma_G\geq1$. Then $I(\Sigma_G)$ is projective. Consequently, every indecomposable injective right $\Sigma_G$-module is projective. This implies that $\Sigma_G$ is
self-injective, contradicting the fact that $G$ is not interval-antipodal. Therefore $\operatorname{domdim}\Sigma_G=0$.
\end{proof}

\subsection{The case of modular lattices}

For notation and conventions concerning modular lattices, we follow Subsection~\ref{subsec:modular-lattice}. We next discuss the homological dimensions of the corresponding distance algebras.

For a finite graph $G$, a \defn{boundary vertex} $y\in V(G)$ is a vertex in $\Max(x)$ for some $x\in V(G)$.
Write $\partial G:=\bigcup_{x\in V(G)}\Max(x)$ for the set of boundary vertices.
We use the standard fact that a finite connected graph with exactly two boundary vertices is a path. 

\begin{thm}\label{thm:modular-lattice-AG}
Let $L$ be a finite modular lattice. The following conditions are
equivalent.
\begin{enumerate}
    \item $\Sigma_L$ is Auslander--Gorenstein.
    \item $\operatorname{domdim}\Sigma_L\geq1$.
    \item $\partial L=Z(L)$.
    \item $L$ is isomorphic to a finite direct product of chains.
    \item $L$ is a singleton or there exist integers $m_1,\ldots,m_c\geq1$ such that
    \[
    L\cong[m_1+1]\times\cdots\times[m_c+1],
    \qquad
    \Sigma_L\cong\bigotimes_{i=1}^c\Sigma_{A_{m_i+1}}.
    \]
\end{enumerate}
When these conditions hold,
\[
\operatorname{idim}\Sigma_L
=\#\{i\mid m_i\geq2\},\; \text{ and }\;
\operatorname{domdim}\Sigma_L=
\begin{cases}
\infty, & m_i=1\text{ for every }i,\\
1, & m_i\geq2\text{ for some }i.
\end{cases}
\]
\end{thm}

\begin{proof}
$(1)\Rightarrow(2)$: This follows by definition. 

$(2)\Rightarrow(3)$:  By Lemma \ref{lem:modular-lattice-central-antipode} (1) $\Leftrightarrow$ (4), we have $Z(L)\subseteq \partial L$.
Now suppose that $y\in\partial L$ and choose an $x\in L$ such that $y\in\Max(x)$. 
Since $\operatorname{soc} P_x$ is given by $S_z$ with $z\in \Max(x)$, the condition $\operatorname{domdim}\Sigma_L\geq1$ means that the direct summand $I_y$ of $I(P_x)$ is projective. 
Hence, it follows from Lemma~\ref{lem:modular-lattice-central-antipode} that $y \in Z(L)$.

$(3)\Rightarrow(4)$:  
Suppose that $L\cong L_1\times L_2\times \cdots\times L_n$ with $L_i$ being directly indecomposable for $1\le i\le n$.  Then the Hasse graph of $L$ is the Cartesian product of those of $L_i$'s .  This implies that 
$\Max_L(x_1,\ldots,x_n)=\prod_{i=1}^n\Max_{L_i}(x_i)$. 
It follows that $\partial L=\partial L_1\times\cdots\times\partial L_n$.
Hence, by Corollary \ref{cor:directly indecomposable lattice}, the assumption (3) implies that $\partial L_i=\{ 0_i < 1_i \}$, which then means that the Hasse diagram of every $L_i$ is a path, and equivalently, every $L_i$ is a chain. 

$(4)\Rightarrow(5)$: The lattice isomorphism is immediate, and the algebra isomorphism follows from Proposition~\ref{prop:tensor}. 

$(5)\Rightarrow(1)$: Each $\Sigma_{A_{m_i+1}}$ is Auslander--Gorenstein by Proposition \ref{prop:idim Sigma tree} and Proposition \ref{prop:domdim Sigma tree}.  Then (1) follows from the fact that tensoring minimal injective coresolutions preserves the Auslander condition.

Finally, the formulae shown follow from the calculations for path graphs and Proposition~\ref{prop:homological-dimensions-tensor}.
\end{proof}

\begin{rmk}\label{rmk:modular-non-IG}
Modularity alone does not imply the Iwanaga--Gorenstein property. 
For any $n\ge 3$, we have a modular lattice whose Hasse diagram is of the form
\begin{center}
\begin{tikzpicture}[
  vertex/.style={
    circle,
    draw,
    fill=white,
    inner sep=0pt,
    minimum size=5mm
  },
  edge/.style={draw=black,line width=.4pt}
]
  \node[vertex] (0) at (0,0) {$0$};
  \node[vertex] (1) at (0,3) {$1$};

  \node[vertex] (x1) at (-3,1.5) {$x_1$};
  \node[vertex] (x2) at (-1.5,1.5) {$x_2$};
  \node (xdots) at (0.5,1.5) {$\cdots\cdots$};
  \node[vertex] (xn) at (3,1.5) {$x_n$};

  \draw[edge] (0) -- (x1);
  \draw[edge] (0) -- (x2);
  \draw[edge] (0) -- (xn);

  \draw[edge] (x1) -- (1);
  \draw[edge] (x2) -- (1);
  \draw[edge] (xn) -- (1);
\end{tikzpicture}
\end{center}
Notice that this is the same as the complete bipartite graph $K_{2,n}$, which is self-centered but not interval-antipodal (see Examples \ref{eg:self-centered eg} and \ref{eg:self-centered-revisited}). 
It follows from Proposition \ref{prop:idim self-centered} that the associated distance algebra has infinite self-injective dimension (so it is not self-injective) and hence is not Iwanaga--Gorenstein; see Example \ref{eg:idim K2,3} for an explicit example.
\end{rmk}

We conclude this subsection with the special case of distributive lattices.  For a finite poset $P$, let $\mathcal J(P)$ denote the lattice of order ideals of $P$, with meet given by intersection and join by union.
Birkhoff's fundamental theorem of finite distributive lattices says that for every finite distributive lattice $L$, we have a lattice isomorphism (poset isomorphism preserving meets and joins) $L\cong \mathcal J(P)$ for some poset $P$.
Moreover, this $P$ is given by the subposet of $L$ consisting of all the join-irreducible elements $x\in L$, i.e., $x\neq 0$ and $x=y\vee z$ implies $x=y$ or $x=z$.

\begin{cor}\label{cor:J(P)-central}
Let $P$ be a finite poset.
\begin{enumerate}
    \item For all $X,Y\in \mathcal J(P)$, $d(X,Y)=|X\mathbin{\triangle}Y|$.
    \item An element $X\in \mathcal J(P)$ is central if and only if it is a union of connected components of the Hasse diagram of $P$.  In such a case, both the central complement and the antipode
of $X$ are $P\setminus X$; in particular, $I(X,P\setminus X)=\mathcal J(P)$.
\end{enumerate}
\end{cor}

\begin{proof}
The rank of an order ideal $X$ is $|X|$. Since meet is intersection,
Lemma~\ref{lem:modular-lattice-distance}~(1) gives
$d(X,Y)=|X|+|Y|-2|X\cap Y|=|X\mathbin{\triangle}Y|$.

Suppose that $X$ is central and let $Y$ be its central complement. Then
$X\cap Y=\emptyset$ and $X\cup Y=P$, so $Y=P\setminus X$. Thus both
$X$ and its complement are order ideals. This is equivalent to saying
that no element of $X$ is comparable with an element of $P\setminus X$,
or, equivalently, that $X$ is a union of connected components of $P$.

Conversely, if $X$ is a union of connected components, then
$P\setminus X$ is also an order ideal, and
\[
\mathcal J(X)\times\mathcal J(P\setminus X)
\longrightarrow\mathcal J(P),
\qquad (U,V)\longmapsto U\cup V,
\]
is a lattice isomorphism. Hence $X$ is central with central complement
$P\setminus X$. The last assertion follows from
Lemma~\ref{lem:modular-lattice-central-antipode}.
\end{proof}

Together with Lemma~\ref{lem:proj-inj-max}, this also shows that the
indecomposable injective $I_X$ is projective if and only if $X$ is a
union of connected components of $P$; in that case,
$I_X\cong P_{P\setminus X}$.

Combining Theorem~\ref{thm:modular-lattice-AG} with Birkhoff's theorem gives the following characterization of the latter property.

\begin{cor}\label{cor:distributive-lattice-AG}
Let $P$ be a finite poset. Then the following conditions are equivalent.
\begin{enumerate}
    \item $\Sigma_{\mathcal{J}(P)}$ is Auslander--Gorenstein.
    \item $\operatorname{domdim}\Sigma_{\mathcal{J}(P)}\geq1$.
    \item Every connected component of $P$ is a chain.
\end{enumerate}
Moreover, $\Sigma_{\mathcal J(P)}$ is minimal Auslander--Gorenstein if and only if every connected component of $P$ is a chain and at most one connected component contains more than one element.
\end{cor}

While the preceding classification of Auslander--Gorenstein
distance algebras admits an elementary proof, classifying
Iwanaga--Gorenstein distance algebras is considerably harder.
For distributive lattices, the main result of the companion
paper \cite{CX2026} yields the following.

\begin{prop}\label{prop:distributive-lattice-IG}
Let $P\neq \emptyset$ be a finite poset. Then the following hold.
\begin{enumerate}[label=\textup{(\arabic*)}]
    \item The distance algebra $\Sigma_{\mathcal J(P)}$ is
    Iwanaga--Gorenstein, with $        \operatorname{idim}\Sigma_{\mathcal J(P)}\leq |P|-1$.
    \item The associated tiled order $\Lambda_{\mathcal J(P)}$
    is Iwanaga--Gorenstein, with
    \[
        \operatorname{idim}\Lambda_{\mathcal J(P)}
        =\operatorname{idim}\Sigma_{\mathcal J(P)}+1
        \leq |P|.
    \]
    \item For any finite modular lattice $L$, the distance algebra $\Sigma_L$ is Iwanaga--Gorenstein if and only if $L$ is distributive.
\end{enumerate}
\end{prop}

\begin{proof}
\textup{(1)} follows directly from
\cite[Theorem~A]{CX2026}. \textup{(2)} follows from \textup{(1)} and Proposition~\ref{prop:homological-dimensions-tiled-order}.

\textup{(3)} Suppose that $\Sigma_L$ is Iwanaga--Gorenstein and $L$ is not distributive. A finite non-distributive modular lattice contains an interval $H=[a,b]$ of length two whose Hasse graph is $K_{2,n}$ for some $n\geq3$. By Lemma~\ref{lem:modular-lattice-distance}, the map $\pi(x)=(x\wedge b)\vee a$ is a projection of $L$ onto $H$. Proposition~\ref{prop:projecting-distance-algebra} implies that $\Sigma_Le_H$ and $e_H\Sigma_L$ are projective as right and left $\Sigma_H$-modules, respectively. Hence the exact functor $(-)e_H$ preserves injective modules. Applying it to a finite injective resolution of $e_H\Sigma_L$ shows that $\Sigma_H\cong e_H\Sigma_Le_H$ has finite right self-injective dimension. The left-module argument is analogous. Thus $\Sigma_H$ is Iwanaga--Gorenstein, contrary to Remark~\ref{rmk:modular-non-IG}.
\end{proof}

\subsection{The case of punctured hypercubes}\label{sec:ex-2-IG}

In the rest of this paper, we show that for any $d\ge 2$, the algebra $\Sigma_{Q_d^-}$ associated with the punctured $d$-dimensional hypercube $Q_d^{-}$ is minimal Auslander--Gorenstein.

Since $Q_d^-$ is a convex subgraph of $Q_d$, the local maxima of the distance function \(d(S,-)\) for $S\subset[d]$ are given by
\[
\Max(S)=
\begin{cases}
\{[d]\setminus S\}, & S\neq\emptyset,\\[2mm]
\bigl\{[d]\setminus\{i\}\mid i\in[d]\bigr\},
& S=\emptyset,
\end{cases}
\]
which means that $P_S \cong I_{S^{\mathrm{c}}} = I_{[d]\setminus S}$ for every non-empty subset $S$ of $[d]$, whereas \(I_{\emptyset}\) is not projective.

\begin{prop}\label{prop:punctured-cube-koszul}
For every $d\geq 1$, the punctured hypercube $Q_d^{-}$ is diagonal.
Consequently, its distance algebra $\Sigma_{Q_d^{-}}$ is Koszul.
\end{prop}

\begin{proof}
We prove the first assertion by induction on $d$. The case $d=1$ is
immediate. Suppose that $d\geq 2$ and that $Q_{d-1}^{-}$ is diagonal.
Write each vertex of $Q_d$ as $(x,i)$, where
$x\in\{0,1\}^{d-1}$ and $i\in\{0,1\}$, and set
\[
Y:=\{(x,0)\mid x\in\{0,1\}^{d-1}\}\cong Q_{d-1},
\]
\[
Z:=\{(x,i)\mid x\neq 1^{d-1},\ i\in\{0,1\}\}
\cong Q_{d-1}^{-}\mathbin{\square}K_2.
\]
Then
\[
Y\cap Z
=\{(x,0)\mid x\neq 1^{d-1}\}
\cong Q_{d-1}^{-},
\qquad
Y\cup Z=Q_d^{-}.
\]
Moreover, $Y\cap Z$ is convex in $Q_d^{-}$, and $Z$ projects onto
$Y\cap Z$ via $\pi(x,i)=(x,0)$. Thus $(Q_d^{-};Y,Z)$ is a projecting decomposition. By the induction hypothesis and Mayer--Vietoris sequences of magnitude (co)homologies \cite[Theorem~6.6]{HW17}, $Q_{d}^{-}$ is diagonal.
\end{proof}

\begin{prop}\label{prop:punctured-distributive-resolution}
Let $P$ be a non-empty finite poset, let $L=\mathcal J(P)$, and regard $L\setminus\{P\}$ as a graph via its undirected Hasse diagram.  Write $\max(P):=\{x\in P\mid \text{there is no $y\in P$ with $x<y$}\}$ for the set of maximal elements of $P$, and set $c:=|\max(P)|$.  Then the indecomposable injective $\Sigma_{L\setminus\{P\}}$-module $I_{\emptyset}$ has a minimal projective resolution
\[
0\longrightarrow P^{-c}\xrightarrow{d_c}P^{-(c-1)}
\longrightarrow\cdots\longrightarrow
P^{-2}\xrightarrow{d_2}P^{-1}
\xrightarrow{\varepsilon}I_{\emptyset}
\longrightarrow0, \text{ where }
P^{-q}:=\bigoplus_{\substack{S\subseteq\max(P)\\|S|=q}} P_{P\setminus S}.
\]
When $c=1$, this means $0\to P^{-1}\xrightarrow{\varepsilon}I_{\emptyset}\to0$.  In particular, $\operatorname{pdim}_{\Sigma_{L\setminus\{P\}}}I_{\emptyset}=c-1$.
\end{prop}

We first explain the connection with \cite[Theorem~2.2]{IM22}.  
Notice that there is a natural algebra quotient $\Sigma_L\to \Bbbk[L^{\operatorname{op}}]=\Bbbk[\mathcal{J}(P^{\operatorname{op}})]$, and then we can identify $I_{\emptyset}$ as the syzygy $N_C=\Omega M_C$ of the antichain module $M_C$ (in the sense of \cite{IM22}) associated to the antichain 
\[C=\{\text{atoms of $L^{\operatorname{op}}$}\}=\{\text{coatoms of $L$}\}=\{P\setminus \{x\}\mid x\in\max(P)\}.\]
The resolution we present here is the `same' as  the one in \cite[Theorem~2.2]{IM22} except that each indecomposable projective $\Bbbk[L^{\operatorname{op}}]$-module is `inflated' to an indecomposable projective $\Sigma_L$-module associated to the same element of $L$.
The proof we present here is basically the same as that of \cite[Theorem~2.2]{IM22}.


\begin{proof}
Fix a total ordering on $\max(P)$.  For $\emptyset\neq S\subseteq\max(P)$ and $x\in S$, the order ideals $P\setminus S$ and $P\setminus(S\setminus\{x\})=(P\setminus S)\sqcup\{x\}$ are adjacent. Let
\[
\lambda_{S,x}:P_{P\setminus S}\longrightarrow
P_{P\setminus(S\setminus\{x\})}
\]
be determined by $\lambda_{S,x}(e_{P\setminus S})=b_{P\setminus(S\setminus\{x\}),\,P\setminus S}$.  For $2\leq q\leq c$, define
$d_q:P^{-q}\to P^{-(q-1)}$
on the summand $P_{P\setminus S}$ by
\[
d_q|_{P_{P\setminus S}}
:=\sum_{x\in S}
(-1)^{|\{y\in S\mid y<x\}|}\lambda_{S,x}.
\]

For $U\in L\setminus\{P\}$, write $\xi_U:=(b_{U,\emptyset})^*\in I_{\emptyset}$.  For $x\in\max(P)$, let
\[
\pi_x:P_{P\setminus\{x\}}\longrightarrow I_{\emptyset}
\]
be determined by $\pi_x(e_{P\setminus\{x\}})=\xi_{P\setminus\{x\}}$, and define $\varepsilon|_{P_{P\setminus\{x\}}}:=\pi_x$.

We verify that these maps form a complex.  Fix $S\subseteq\max(P)$ with $|S|\geq3$ and distinct $x,y\in S$, say $x<y$.  The component of $d_{q-1}d_q$ from $P_{P\setminus S}$ to $P_{P\setminus(S\setminus\{x,y\})}$ is the sum of the two composites obtained by deleting $x$ first or $y$ first.  Both send
$e_{P\setminus S}$ to $b_{P\setminus(S\setminus\{x,y\}),\,P\setminus S}$, since the two length-two paths are shortest paths with the same endpoints.  Put $n_S(z):=|\{w\in S\mid w<z\}|$.  The signs of the two composites are
$(-1)^{n_S(x)+n_S(y)-1}$ and $(-1)^{n_S(x)+n_S(y)}$, respectively, and hence cancel.  Thus $d_{q-1}d_q=0$.

If $S=\{x<y\}$, the two components of $d_2|_{P_{P\setminus S}}$ have coefficients $1$ and $-1$, while their composites with $\varepsilon$ both send $e_{P\setminus S}$ to $\xi_{P\setminus S}$.  Hence $\varepsilon d_2=0$, and the displayed sequence is a complex.

We prove exactness after applying $(-)\cdot e_U$ for a fixed $U\in L\setminus\{P\}$.  Set $R_U:=\max(P)\setminus U$ and $T_U:=\max(P)\cap U$.  Since $U$ is a
proper order ideal, $R_U$ is non-empty.  Moreover, $(P_{P\setminus S})e_U=\Bbbk b_{P\setminus S,U}$ and $I_{\emptyset}e_U=\Bbbk\xi_U$.

For $x\in S$, multiplication in the distance algebra gives
\[
\lambda_{S,x}(b_{P\setminus S,U})
=\begin{cases}
b_{P\setminus(S\setminus\{x\}),U},&x\in R_U,\\
0,&x\in T_U.
\end{cases}
\]
Indeed, adjoining $x$ to $P\setminus S$ increases its distance from
$U$ by one exactly when $x\notin U$.  Similarly,
\[
\pi_x(b_{P\setminus\{x\},U})
=\begin{cases}
\xi_U,&x\in R_U,\\
0,&x\in T_U.
\end{cases}
\]

Fix $F\subseteq T_U$ and collect the basis elements indexed by $S=F\sqcup B$ with $B\subseteq R_U$.  The preceding formulas show that the differential preserves $F$ and only deletes elements of $B$.  If $F\neq\emptyset$, the vector corresponding to $B=\emptyset$ is $b_{P\setminus F,U}$; if $F=\emptyset$, it is $\xi_U$.

For $x\in R_U$, set $a_F(x):=|\{y\in F\mid y<x\}|$.  For $B\subseteq R_U$, define
\[
\eta_{F,B}:=
\begin{cases}
(-1)^{\sum_{x\in B}a_F(x)}
b_{P\setminus(F\sqcup B),U},&F\sqcup B\neq\emptyset,\\
\xi_U,&F=B=\emptyset.
\end{cases}
\]
The definition of the differentials now gives
\[
d(\eta_{F,B})
=\sum_{x\in B}
(-1)^{|\{y\in B\mid y<x\}|}
\eta_{F,B\setminus\{x\}}.
\]

Choose the least element $r$ of $R_U$ and define
\[
h(\eta_{F,B})
:=
\begin{cases}
\eta_{F,B\sqcup\{r\}},&r\notin B,\\
0,&r\in B.
\end{cases}
\]
If $r\notin B$, the term obtained by deleting $r$ from $dh(\eta_{F,B})$ is $\eta_{F,B}$, and the remaining terms cancel with $hd(\eta_{F,B})$.  If $r\in B$, only the term obtained by deleting $r$ survives after applying $h$, and it is sent to $\eta_{F,B}$.  Thus $dh+hd=\operatorname{id}$.  Each summand indexed by $F$ is contractible, so the displayed complex is exact after
applying $-e_U$.  Since this holds for every vertex $U$, the complex
is exact.

Finally, every non-zero entry of $d_q$ is induced by an arrow, so $\operatorname{Im}d_q\subseteq\operatorname{rad}P^{-(q-1)}$.  Exactness gives $\ker\varepsilon=\operatorname{Im}d_2\subseteq\operatorname{rad}P^{-1}$ when $c\geq2$; for $c=1$, the kernel is zero.  Hence $\varepsilon$ is a projective cover, and the displayed resolution is minimal.  Since $P^{-c}=P_{P\setminus\max(P)}\neq0$, its length is $c-1$.
\end{proof}

\begin{rmk}
One can ask if the same description holds for the minimal projective resolution of $I_{0}$ when we replace the distributive lattice $L$ by a finite modular lattice.  The answer is no.
For example, if $L$ is the five-element modular lattice with atoms $a,b,c$ (the Hasse graph being isomorphic to $K_{2,3}$), then, after removing its maximum $\hat 1$, one obtains
$$0\longrightarrow P_{\hat 0}^{\oplus 2}
\longrightarrow P_a\oplus P_b\oplus P_c
\longrightarrow I_{\hat 0}
\longrightarrow 0.$$
\end{rmk}

We thus obtain the following result.

\begin{cor}\label{prop:punctured-cube-injective-resolution}
Let $d\geq 2$. Then the unique non-projective indecomposable injective module $I_{\emptyset}$ has a minimal projective resolution
$$0\longrightarrow P_{\emptyset}
\longrightarrow
\bigoplus_{\substack{S\subsetneq[d]\\ |S|=1}}P_S
\longrightarrow\cdots\longrightarrow
\bigoplus_{\substack{S\subsetneq[d]\\ |S|=d-1}}P_S
\longrightarrow I_{\emptyset}
\longrightarrow0.$$
Dually, this gives a minimal injective coresolution of $P_{\emptyset}$. Therefore,
$$\operatorname{findim}\Sigma_{Q_d^{-}}=\operatorname{domdim}\Sigma_{Q_d^{-}}=\operatorname{idim}\Sigma_{Q_d^{-}}=d-1.$$
Thus $\Sigma_{Q_d^{-}}$ is $(d-2)$-minimal Auslander--Gorenstein.
\end{cor}

\begin{eg}\label{ex:hy-cube-d=2}
The graph $Q_3^{-}$ is as follows:
$$\begin{tikzpicture}[
scale=0.6,
    every node/.style={
        circle,
        draw,
        minimum size=8mm,
        inner sep=0pt
    }
]
\node (empty) at (0,0)       {$\emptyset$};
\node (one)   at (0,2.6)     {$\{1\}$};
\node (onetwo) at (2.25,1.3) {$\{1,2\}$};
\node (two)   at (2.25,-1.3) {$\{2\}$};
\node (twothree) at (0,-2.6) {$\{2,3\}$};
\node (three) at (-2.25,-1.3) {$\{3\}$};
\node (onethree) at (-2.25,1.3) {$\{1,3\}$};

\draw (one) -- (onetwo);
\draw (onetwo) -- (two);
\draw (two) -- (twothree);
\draw (twothree) -- (three);
\draw (three) -- (onethree);
\draw (onethree) -- (one);

\draw (empty) -- (one);
\draw (empty) -- (two);
\draw (empty) -- (three);
\end{tikzpicture}$$
A direct calculation shows that there is a minimal injective resolution
\begin{align*}
 0\longrightarrow  P_{\emptyset} & \longrightarrow P^{-1} \longrightarrow P^0 \longrightarrow I_{\emptyset} \longrightarrow 0,\\
\text{ where }  P^{-1} & = \bigoplus_{|X|=1} P_X = P_{\{1\}}\oplus P_{\{2\}}\oplus P_{\{3\}} \cong  I_{\{2,3\}}\oplus I_{\{1,3\}}\oplus I_{\{1,2\}}, \\
\text{ and } {P^0} & = \bigoplus_{|X|=2} P_X = P_{\{2,3\}}\oplus P_{\{1,3\}}\oplus P_{\{1,2\}} \cong I_{\{1\}}\oplus I_{\{2\}}\oplus I_{\{3\}}. 
\end{align*}
Hence, we have
$\operatorname{findim}\Sigma_{Q_3^{-}}
=
\operatorname{domdim}\Sigma_{Q_3^{-}}
=
\operatorname{idim}\Sigma_{Q_3^{-}}
=
2.$
\end{eg}

\subsubsection{The associated precluster tilting module}\label{subsubsec:precluster}

Recall that, for $n\geq1$, an \defn{$n$-precluster tilting module} over an algebra $\Lambda$ is a generator-cogenerator $M$ satisfying
$\operatorname{Ext}_{\Lambda}^i(M,M)=0$ for $0<i<n$ and
$\tau_nM,\tau_n^{-}M\in\operatorname{add}M$, where
$\tau_n:=\tau\Omega^{n-1}$ and
$\tau_n^{-}:=\tau^{-}\Omega^{-(n-1)}$ are the higher
Auslander--Reiten translations; see \cite[Definition~3.2]{IS18}.

The \defn{Auslander--Iyama--Solberg correspondence} identifies
$n$-minimal Auslander--Gorenstein algebras with endomorphism
algebras of $n$-precluster tilting modules. More precisely, in our
right-module convention, if $\Gamma$ is $n$-minimal
Auslander--Gorenstein and $I$ is an additive generator of its
projective-injective modules, then for
$\Lambda:=\operatorname{End}_{\Gamma}(I)$,
the module $M:=\operatorname{Hom}_{\Gamma}(I,\Gamma)$ is an $n$-precluster
tilting $\Lambda$-module satisfying
$\Gamma\cong\operatorname{End}_{\Lambda}(M)$; see
\cite[Proposition~4.4 and Theorem~4.5]{IS18}.

We apply this correspondence to $\Sigma_{Q_d^{-}}$ for $d\geq3$. Put
$f:=1-e_{\emptyset}$ and $B:=f\Sigma_{Q_d^{-}}f$. 
By Proposition~\ref{prop:punctured-cube-injective-resolution}, $\Sigma_{Q_d^{-}}$ is
$(d-2)$-minimal Auslander--Gorenstein, and $D(\Sigma_{Q_d^-}f)$ is an additive
generator of its projective-injective right $\Sigma_{Q_d^{-}}$-modules.

\begin{prop}\label{prop:punctured-cube-precluster}
The following statements hold.
\begin{enumerate}
    \item $B$ is isomorphic, as an algebra, to the self-injective algebra $\Sigma_{Q_d^{
\pm}}$, where $Q_d^{\pm}=Q_d\setminus\{\emptyset,[d]\}$ is the    doubly punctured hypercube. 

    \item Let $N:=e_{\emptyset}\Sigma_{Q_d^{-}}f=\operatorname{rad}_{\Sigma_{Q_d^{-}}} P_{\emptyset}$.
    Then for every $d>2$, the right $B$-module $M:=\Sigma_{Q_d^{-}}f\cong B\oplus N$ is a basic $(d-2)$-precluster
    tilting $B$-module.  In particular, we have $\tau_{d-2}N\cong N\cong\tau_{d-2}^{-}N$ and $\Sigma_{Q_d^{-}}\cong\operatorname{End}_B(M)$ is a Morita algebra in the sense of \cite{KY13}.  
\end{enumerate}
\end{prop}

\begin{proof}
(1) Recall from Example~\ref{ex:interval-antipodals}(2) that $Q_d^{\pm}$ is an interval-antipodal convex full subgraph of $Q_d^{-}$. Hence
Lemma~\ref{lem:corner-convex-subgraph} gives
$B\cong\Sigma_{Q_d^{\pm}}$, and self-injectivity is just Theorem~\ref{thm:self-inj}.

(2) The non-projective direct summand of the $B$-module $\Sigma_{Q_d^-}f$ is given by $(1-f)\Sigma_{Q_d^-} f = N$.  
Since $\Sigma_{Q_d^-}$ is minimal Auslander--Gorenstein, the rest of the statement follows immediately from the theory of the Auslander--Iyama--Solberg correspondence \cite[Proposition~4.4, Theorem~4.5, Section 6]{IS18}.
\end{proof}

\bibliographystyle{alpha}
\bibliography{reference.bib}

@article{AHHT,
  author  = {Angeleri H{\"u}gel, L. and Herbera, D. and Trlifaj, J.},
  title   = {Tilting modules and {G}orenstein rings},
  journal = {Forum Math.},
  volume  = {18},
  number  = {2},
  pages   = {211--229},
  year    = {2006},
  doi     = {10.1515/FORUM.2006.013}
}

@misc{AI24,
  author        = {Asao, Y. and Ivanov, S. O.},
  title         = {Magnitude homology is a derived functor},
  year          = {2024},
  eprint        = {2402.14466},
  archivePrefix = {arXiv},
  note          = {arXiv:2402.14466}
}

@book{ARS95,
  author    = {Auslander, M. and Reiten, I. and Smal{\o}, S. O.},
  title     = {Representation theory of {A}rtin algebras},
  series    = {Camb. Stud. Adv. Math.},
  volume    = {36},
  publisher = {Cambridge Univ. Press},
  address   = {Cambridge},
  year      = {1995}
}

@article{AW24,
  author  = {Asao, Y. and Wakatsuki, S.},
  title   = {Minimal projective resolution and magnitude homology of geodetic metric spaces},
  journal = {Homology Homotopy Appl.},
  volume  = {28},
  number  = {2},
  pages   = {123--154},
  year    = {2026},
  doi     = {10.4310/hha.2026.v28.n2.a7}
}

@book{BH90,
  author    = {Buckley, F. and Harary, F.},
  title     = {Distance in graphs},
  publisher = {Addison-Wesley},
  address   = {Redwood City, CA},
  year      = {1990}
}

@article{BK21,
  author  = {Bottinelli, R. and Kaiser, T.},
  title   = {Magnitude homology, diagonality, and median spaces},
  journal = {Homology Homotopy Appl.},
  volume  = {23},
  number  = {2},
  pages   = {121--140},
  year    = {2021},
  doi     = {10.4310/hha.2021.v23.n2.a7}
}

@book{BM08,
  author    = {Bondy, J. A. and Murty, U. S. R.},
  title     = {Graph theory},
  series    = {Grad. Texts in Math.},
  volume    = {244},
  publisher = {Springer},
  address   = {New York},
  year      = {2008}
}

@article{BM10,
  author  = {Bergh, P. A. and Madsen, D.},
  title   = {Hochschild homology and split pairs},
  journal = {Bull. Sci. Math.},
  volume  = {134},
  number  = {7},
  pages   = {665--676},
  year    = {2010}
}

@article{BR87,
  author  = {Butler, M. C. R. and Ringel, C. M.},
  title   = {{A}uslander--{R}eiten sequences with few middle terms and applications to string algebras},
  journal = {Comm. Algebra},
  volume  = {15},
  number  = {1--2},
  pages   = {145--179},
  year    = {1987},
  doi     = {10.1080/00927878708823416}
}

@article{Buc89,
  author  = {Buckley, F.},
  title   = {Self-centered graphs},
  journal = {Ann. N.Y. Acad. Sci.},
  volume  = {576},
  number  = {1},
  pages   = {71--78},
  year    = {1989},
  doi     = {10.1111/j.1749-6632.1989.tb16384.x}
}

@article{CB95,
  author  = {Crawley-Boevey, W.},
  title   = {Tameness of biserial algebras},
  journal = {Arch. Math. (Basel)},
  volume  = {65},
  number  = {5},
  pages   = {399--407},
  year    = {1995},
  doi     = {10.1007/BF01198070}
}

@book{CE56,
  author    = {Cartan, H. and Eilenberg, S.},
  title     = {Homological algebra},
  publisher = {Princeton Univ. Press},
  address   = {Princeton, NJ},
  year      = {1956}
}

@unpublished{CX2026,
  author = {Chan, A. and Xing, B.},
  title  = {Distance algebras of distributive lattices are {G}orenstein},
  note   = {In preparation},
  year   = {2026}
}

@article{DKR67,
  author  = {Drozd, Yu. A. and Kirichenko, V. V. and Roiter, A. V.},
  title   = {On hereditary and {B}ass orders},
  journal = {Math. USSR-Izv.},
  volume  = {1},
  number  = {6},
  pages   = {1357--1375},
  year    = {1967},
  doi     = {10.1070/IM1967v001n06ABEH000625}
}

@misc{EKM26,
  author        = {Easwar, R. and Kuber, A. and Mittal, M.},
  title         = {On zeta functions and $\mu$-series of string algebras},
  year          = {2026},
  eprint        = {2604.03108},
  archivePrefix = {arXiv},
  primaryClass  = {math.RT},
  note          = {arXiv:2604.03108}
}

@article{ERZ57,
  author  = {Eilenberg, S. and Rosenberg, A. and Zelinsky, D.},
  title   = {On the dimension of modules and algebras, {VIII}: dimension of tensor products},
  journal = {Nagoya Math. J.},
  volume  = {12},
  pages   = {71--93},
  year    = {1957},
  doi     = {10.1017/S0027763000021954}
}

@article{ES92,
  author  = {Erdmann, K. and Skowro{\'n}ski, A.},
  title   = {On {A}uslander--{R}eiten components of blocks and self-injective biserial algebras},
  journal = {Trans. Amer. Math. Soc.},
  volume  = {330},
  number  = {1},
  pages   = {165--189},
  year    = {1992}
}

@article{Fuller68,
  author  = {Fuller, K. R.},
  title   = {Biserial rings},
  journal = {Lecture Notes in Math.},
  volume  = {87},
  pages   = {64--90},
  year    = {1969}
}

@article{GP79,
  author  = {Gelfand, I. M. and Ponomarev, V. A.},
  title   = {Model algebras and representations of graphs},
  journal = {Funct. Anal. Appl.},
  volume  = {13},
  number  = {3},
  pages   = {157--166},
  year    = {1979}
}

@incollection{GHZ98,
  author    = {Goodearl, K. R. and Huisgen-Zimmermann, B.},
  title     = {Repetitive resolutions over classical orders and finite-dimensional algebras},
  booktitle = {Algebras and modules, {II}},
  series    = {CMS Conf. Proc.},
  volume    = {24},
  pages     = {205--225},
  publisher = {Amer. Math. Soc.},
  address   = {Providence, RI},
  year      = {1998}
}

@book{Gra03,
  author    = {Gr{\"a}tzer, G.},
  title     = {General lattice theory},
  edition   = {2},
  publisher = {Birkh{\"a}user},
  address   = {Basel},
  year      = {2003}
}

@article{GS16,
  author  = {Green, E. L. and Schroll, S.},
  title   = {Multiserial and special multiserial algebras and their representations},
  journal = {Adv. Math.},
  volume  = {302},
  pages   = {1111--1136},
  year    = {2016}
}

@article{Han06,
  author  = {Han, Y.},
  title   = {Hochschild {(co)homology} dimension},
  journal = {J. Lond. Math. Soc. (2)},
  volume  = {73},
  number  = {3},
  pages   = {657--668},
  year    = {2006}
}

@article{Has51,
  author  = {Hashimoto, J.},
  title   = {On direct product decomposition of partially ordered sets},
  journal = {Ann. of Math. (2)},
  volume  = {54},
  pages   = {315--318},
  year    = {1951},
  doi     = {10.2307/1969532}
}

@article{Her08,
  author  = {Herschend, M.},
  title   = {Tensor products on quiver representations},
  journal = {J. Pure Appl. Algebra},
  volume  = {212},
  number  = {2},
  pages   = {452--469},
  year    = {2008}
}

@book{HIK11,
  author    = {Hammack, R. and Imrich, W. and Klav{\v{z}}ar, S.},
  title     = {Handbook of product graphs},
  edition   = {2},
  publisher = {CRC Press},
  address   = {Boca Raton, FL},
  year      = {2011}
}

@article{HW17,
  author  = {Hepworth, R. and Willerton, S.},
  title   = {Categorifying the magnitude of a graph},
  journal = {Homology Homotopy Appl.},
  volume  = {19},
  number  = {2},
  pages   = {31--60},
  year    = {2017},
  doi     = {10.4310/HHA.2017.v19.n2.a3}
}

@article{HZ91,
  author  = {Huisgen-Zimmermann, B.},
  title   = {Predicting syzygies over monomial relation algebras},
  journal = {Manuscripta Math.},
  volume  = {70},
  number  = {2},
  pages   = {157--182},
  year    = {1991}
}

@misc{IL26,
  author        = {Iyama, O. and Liu, J.},
  title         = {Triangle equivalences between {G}orenstein tiled orders and incidence algebras of posets},
  year          = {2026},
  eprint        = {2602.01878},
  archivePrefix = {arXiv},
  primaryClass  = {math.RT},
  note          = {arXiv:2602.01878}
}

@article{IM22,
  author  = {Iyama, O. and Marczinzik, R.},
  title   = {Distributive lattices and {A}uslander regular algebras},
  journal = {Adv. Math.},
  volume  = {398},
  pages   = {108233},
  year    = {2022},
  doi     = {10.1016/j.aim.2022.108233}
}

@article{IM24,
  author  = {Ivanov, S. O. and Mukoseev, L.},
  title   = {On diagonal digraphs, {K}oszul algebras and triangulations of homology spheres},
  journal = {Adv. Math.},
  volume  = {502},
  pages   = {111163},
  year    = {2026},
  note    = {Part B},
  doi     = {10.1016/j.aim.2026.111163}
}

@article{IS18,
  author  = {Iyama, O. and Solberg, {\O}.},
  title   = {Auslander--{G}orenstein algebras and precluster tilting},
  journal = {Adv. Math.},
  volume  = {326},
  pages   = {200--240},
  year    = {2018},
  doi     = {10.1016/j.aim.2017.11.025}
}

@article{Iya07,
  author  = {Iyama, O.},
  title   = {Auslander correspondence},
  journal = {Adv. Math.},
  volume  = {210},
  number  = {1},
  pages   = {51--82},
  year    = {2007},
  doi     = {10.1016/j.aim.2006.06.003}
}

@article{KY13,
  author  = {Kerner, O. and Yamagata, K.},
  title   = {Morita algebras},
  journal = {J. Algebra},
  volume  = {382},
  pages   = {185--202},
  year    = {2013},
  doi     = {10.1016/j.jalgebra.2013.02.013}
}

@book{Lod98,
  author    = {Loday, J.-L.},
  title     = {Cyclic homology},
  series    = {Grundlehren Math. Wiss.},
  volume    = {301},
  edition   = {2},
  publisher = {Springer},
  address   = {Berlin},
  year      = {1998}
}

@article{Miz14,
  author  = {Mizuno, Y.},
  title   = {Classifying {$\tau$}-tilting modules over preprojective algebras of {D}ynkin type},
  journal = {Math. Z.},
  volume  = {277},
  number  = {3},
  pages   = {665--690},
  year    = {2014}
}

@article{Mul68,
  author  = {M{\"u}ller, B.},
  title   = {The classification of algebras by dominant dimension},
  journal = {Canad. J. Math.},
  volume  = {20},
  pages   = {398--409},
  year    = {1968}
}

@article{MV99,
  author  = {Mart{\'i}nez-Villa, R.},
  title   = {Graded, selfinjective, and {K}oszul algebras},
  journal = {J. Algebra},
  volume  = {215},
  number  = {1},
  pages   = {34--72},
  year    = {1999},
  doi     = {10.1006/jabr.1998.7728}
}

@article{Rog73,
  author  = {Roggenkamp, K. W.},
  title   = {Injective modules for group rings and {G}orenstein orders},
  journal = {J. Algebra},
  volume  = {24},
  number  = {3},
  pages   = {465--472},
  year    = {1973},
  doi     = {10.1016/0021-8693(73)90120-8}
}

@incollection{Smi96,
  author    = {Smith, S. P.},
  title     = {Some finite-dimensional algebras related to elliptic curves},
  booktitle = {Representation theory of algebras and related topics},
  series    = {CMS Conf. Proc.},
  volume    = {19},
  pages     = {315--348},
  publisher = {Amer. Math. Soc.},
  address   = {Providence, RI},
  year      = {1996}
}

@article{SW83,
  author  = {Skowro{\'n}ski, A. and Waschb{\"u}sch, J.},
  title   = {Representation-finite biserial algebras},
  journal = {J. Reine Angew. Math.},
  volume  = {345},
  pages   = {172--181},
  year    = {1983}
}

@book{Wei94,
  author    = {Weibel, C. A.},
  title     = {An introduction to homological algebra},
  series    = {Camb. Stud. Adv. Math.},
  volume    = {38},
  publisher = {Cambridge Univ. Press},
  address   = {Cambridge},
  year      = {1994},
  doi       = {10.1017/CBO9781139644136}
}

@article{Wil83,
  author  = {Wilson, G. V.},
  title   = {The {C}artan map on categories of graded modules},
  journal = {J. Algebra},
  volume  = {85},
  number  = {2},
  pages   = {390--398},
  year    = {1983}
}

@misc{Xin24,
  author        = {Xing, B.},
  title         = {Quasi-biserial algebras, special quasi-biserial algebras and symmetric fractional {B}rauer graph algebras},
  year          = {2024},
  eprint        = {2408.03778},
  archivePrefix = {arXiv},
  note          = {arXiv:2408.03778}
}
\Addresses

\end{document}